\documentclass[12pt]{amsart}
\usepackage{amscd}
\usepackage{pinlabel}
\usepackage{graphicx}

\newtheorem{thm}{Theorem}[section]

\newtheorem{dfn}[thm]{Definition}

\newtheorem{cor}[thm]{Corollary}

\newtheorem{lma}[thm]{Lemma}

\newtheorem{rem}[thm]{Remark}
\newenvironment{rmk}{\begin{rem}\rm}{\end{rem}}

\newenvironment{pf}{\begin{proof}}{\end{proof}}

\numberwithin{equation}{section}

\newcommand{\R}{{\mathbb{R}}}
\newcommand{\C}{{\mathbb{C}}}

\newcommand{\Z}{{\mathbb{Z}}}
\newcommand{\CC}{{\mathcal{C}}}

\newcommand{\A}{{\mathcal{A}}}
\newcommand{\Pe}{{\mathcal{P}}}
\newcommand{\Ordo}{{\mathbf{O}}}

\newcommand{\M}{{\mathcal{M}}}
\newcommand{\Co}{{\mathcal{G}}}
\newcommand{\Rch}{{\mathcal{R}}}
\newcommand{\Fu}{{\mathcal{F}}}

\newcommand{\ev}{\operatorname{ev}}
\newcommand{\la}{\langle}
\newcommand{\ra}{\rangle}
\newcommand{\pa}{\partial}
\newcommand{\id}{\operatorname{id}}

\newcommand{\krn}{\operatorname{ker}}

\newcommand{\img}{\operatorname{Im}}

\title[Isotopies of Legendrian 1-knots and Legendrian 2-tori]
{Isotopies of Legendrian 1-knots\\ and Legendrian 2-tori}
\author{Tobias Ekholm}
\author{Tam{\'a}s K{\'a}lm{\'a}n}
\address{Department of mathematics, Uppsala University, Box 480, 751 06 Uppsala, Sweden}
\email{tobias{\@@}math.uu.se}
\address{USC, Department of mathematics, 3620 S Vermont Ave, Los Angeles, CA 90089}
\email{tkalman{\@@}usc.edu}

\begin{document}

\begin{abstract}
We construct a Legendrian $2$-torus in the $1$-jet space of $S^1\times\R$ (or of $\R^2$) from a loop of Legendrian knots in the $1$-jet space of $\R$. The differential graded algebra (DGA) for the Legendrian contact homology of the torus is explicitly computed in terms of the DGA of the knot and the monodromy operator of the loop. The contact homology of the torus is shown to depend only on the chain homotopy type of the monodromy operator. The construction leads to many new examples of Legendrian knotted tori. In particular, it allows us to construct a Legendrian torus with DGA which does not admit any augmentation (linearization) but which still has non-trivial homology, as well as two Legendrian tori with isomorphic linearized contact homologies but with distinct contact homologies.
\end{abstract}

\thanks{TE is an Alfred P. Sloan Research Fellow, acknowledges support from NSF-grant DMS-0505076, and from the Swedish Royal Academy of Sciences, the Knut and Alice Wallenberg foundation.}
\maketitle
%\tableofcontents

\section{Introduction}\label{S:intr}
Let $M$ be a smooth $n$-manifold and consider the $1$--jet space $J^1(M)=T^\ast M\times\R$ with coordinates $(x,y,z)$, where $x\in M$, $y\in T_x^\ast M$, and $z\in\R$. The $1$--form $\alpha=dz-\sum_j y_j\,dx_j$ is a contact $1$--form and $\xi=\krn(\alpha)$ is the standard contact structure on $J^1(M)$. The Reeb vector field of a contact form $\beta$ is the unique vector field $R_\beta$ which satisfies $d\beta(R_\beta,\cdot)=0$ and $\beta(R_\beta)=0$. Hence, if $\alpha$ is the standard contact form then $R_\alpha=\pa_z$. An $n$--dimensional submanifold $L\subset J^1(M)$ which is everywhere tangent to $\xi$ is called {\em Legendrian}. The isotopy problem for Legendrian submanifolds is the problem of distinguishing the path components of the space of Legendrian submanifolds. It is a very rich and interesting problem, see e.g. \cite{Ch, El, EES2, Ng1, Ng2}. In \cite{K} the related problem of detecting non-contractible loops in the space of Legendrian $1$-submanifolds of $J^1(\R)$ was studied. In this paper we relate the latter problem to the former problem for Legendrian tori in $J^1(\R\times S^1)$ (and in $J^{1}(\R^{2})$).

The introduction of Legendrian contact homology, see \cite{Ch, EGH}, lead to major breakthroughs in the study of the Legendrian isotopy problem. The contact homology of a Legendrian submanifold $L\subset J^{1}(M)$ is the homology of a differential graded algebra (DGA) $\A(L)$ associated to $L$. Here $\A(L)$ is the (non-commutative) unital DGA freely generated by the Reeb chords of $L$ (flow lines of the Reeb vector field starting and ending on $L$) graded by their Conley-Zhender indices. (For simplicity, throughout this paper we will work with $\Z_2$-coefficients so that $\A(L)$ is an algebra over a $\Z_2$-group ring, rather than a $\Z$-group ring.) The differential on $\A(L)$ is defined using moduli spaces of boundary punctured holomorphic disks in the symplectization $J^1(M)\times\R$ with Lagrangian boundary condition $L\times\R\subset J^1(M)\times\R$ and with certain asymptotic properties near the punctures, see Subsections \ref{s:hodi} and \ref{s:LCH}. For $1$--dimensional Legendrian submanifolds of $J^1(\R)$, the Riemann mapping theorem allows for a purely combinatorial description of contact homology, see \cite{Ch}. For higher dimensional Legendrian submanifolds, analytical foundations for contact homology were worked out in \cite{EES2, EES4}. In \cite{E}, a more combinatorial description was obtained: Legendrian contact homology in $1$--jet spaces was described entirely in terms of Morse theoretic objects called flow trees, see Subsection \ref{s:fltr}.

Contact homology has certain functorial properties. In particular, each isotopy of Legendrian knots induces a Lagrangian cobordism which in turn induces a morphism of the contact homologies at its endpoints. Such morphisms induced by isotopies were dealt with from a purely combinatorial point of view in \cite{K}. (The techniques used in this paper lead to a description of these morphisms in terms of moduli spaces of holomorphic curves in a cobordism, see \cite{EHK} ).

In order to state the main theorem of the paper we first describe the underlying geometric construction. Let $L\subset J^1(M)$ be a Legendrian submanifold which is in sufficiently general position with respect to the projection $\Pi_F\colon J^1(M)\to J^0(M)$. The image $\Pi_F(L)$ is the {\em front of $L$} and determines $L$. If $\gamma(t)$, $t\in S^1$ is a $1$--parameter family of Legendrian submanifolds, starting and ending at $L$, then $(\Pi_F(\gamma(t)),t)\subset J^0(M\times S^1)$ is the front of a Legendrian embedding $L\times S^1\to J^1(M\times S^1)$. We call this Legendrian embedding the {\em trace of the isotopy} and denote the corresponding Legendrian submanifold by $\Sigma_\gamma(L)$. When $M=\R$, consider an embedding $S^{1}\times\R\to\R^{2}$ obtained by identifying $S^{1}\times\R$ with a small tubular neighborhood of the unit circle. Using this embedding we may consider $\Sigma_\gamma(L)\subset J^{1}(S^{1}\times\R)\subset J^{1}(\R^{2})$.

\begin{thm}\label{t:main}
Let $\gamma(t)$, $t\in S^1$ be a loop of Legendrian submanifolds of $J^{1}(\R)$ starting and ending at $L\subset J^1(\R)$. Let $\A=\Z_{2}[H_1(L)]\la c_1,\dots,c_r\ra$ be the DGA of $L$, let $|c_j|$ denote the degree of the generator $c_j$, and let $\phi\colon\A\to\A$ be the monodromy operator associated to $\gamma$. Then the DGA of $\Sigma_\gamma(L)$ as a Legendrian submanifold of $J^{1}(S^{1}\times\R)$ or of $J^{1}(\R^{2})$ is stable tame isomorphic to the algebra $(\hat\A,\Delta)$, where
\[
\hat\A=\Z_2[H_1(L\times S^1)]\la c_1,\dots,c_r,\hat c_1,\dots,\hat c_r\ra,\\
\]
with $|\hat c_j|=|c_j|+1$ for all $j$, and where
\begin{align*}
\Delta(c_j)&=\pa c_j,\\
\Delta(\hat c_j)&= c_j+\phi(c_j)+\Gamma_\phi(\pa c_j),
\end{align*}
with $\Gamma_\phi\colon \A\to\hat \A$ denoting the degree $1$ derivation defined by
$\Gamma(c_s)=\hat c_s$ and extended to all of $\A$ by
$$
\Gamma_\phi(\alpha\beta)=\Gamma_\phi(\alpha)\phi(\beta)+\alpha\Gamma_\phi(\beta).
$$
\end{thm}
Theorem \ref{t:main} is proved in Section \ref{S:algebra}. Using it, we produce examples establishing the following.
\begin{thm}\label{t:ex}
Let $Y$ denote $J^{1}(S^{1}\times\R)$ or $J^{1}(\R^{2})$ equipped with its standard contact form.
\begin{itemize}
\item[{\rm (a)}]
There exists a Legendrian torus in $Y$ with DGA which has no augmentation but which has non-zero homology.
\item[{\rm (b)}]
There exists a Legendrian torus in $Y$ with linearized contact homology isomorphic to the linearized contact homology of the standard torus (the standard torus is the trace of the constant isotopy of the unknot) but with full contact homology different from that of the standard torus.
\end{itemize}
\end{thm}
Theorem \ref{t:ex} is proved in Section \ref{S:examples}.

The paper is organized as follows. In Section \ref{S:background} we recall the definition of contact homology and describe how to compute it in $1$--jet spaces using flow trees. In Section \ref{S:gench} we discuss certain slight generalizations of contact homology described in terms of perturbed flow trees. In Section \ref{S:geomprt} we explain how the differential $\pa$ of the DGA $\A$ of the trace of the constant isotopy can be computed in terms of perturbed flow trees which arise from a geometric perturbation of the Legendrian submanifold. This computation is however not explicit enough to yield a closed formula for $\pa$. In Section \ref{S:abstrprt} we design an abstract perturbation for the trace of the constant isotopy which yields a closed formula for another differential $\pa'$ on the DGA $\A$. By a result from Section \ref{S:gench}, the DGAs $(\A,\pa)$ and $(\A,\pa')$ are tame isomorphic. In Section \ref{S:move} we discuss how to decompose an isotopy into simple pieces and how to compute the DGAs of the traces of such simple pieces. Concatenation of the pieces then yields the DGA for the full isotopy. However, the concatenation gives a DGA with an enormous number of generators. In Section \ref{S:algebra} we show, using a purely algebraic argument, that the DGA with the large number of generators is stable tame isomorphic to the DGA in the formulation of Theorem \ref{t:main}. We also give an algebraic proof of fact that the contact homology of the trace of an isotopy depends only on the chain homotopy class of the morphism of the isotopy. In Section \ref{S:examples} we study examples needed to establish Theorem \ref{t:ex}. 
\section{Background}\label{S:background}
In this section we give a brief description of Legendrian contact homology and of how to compute it for Legendrian submanifolds in $1$--jet spaces. For details we refer to \cite{EES2, EES4, E}. More precisely, the DGA associated to a Legendrian submanifold. is described in Subsections \ref{s:hodi} and \ref{s:LCH}. In Subsection \ref{s:fltr} we define flow trees and describe the relation between rigid flow trees and rigid holomorphic disks.

\subsection{Holomorphic disks}\label{s:hodi}
Let $M$ be a smooth $n$-manifold, endow $J^1(M)=T^\ast M\times\R$ with its standard contact form, and let $z$ be a coordinate in the $\R$-direction. Let $L$ be a closed Legendrian submanifold . We assume that $L$ is sufficiently generic so that the {\em Lagrangian projection} $\Pi_\C\colon J^1(M)\to T^\ast M$ restricted to $L$ has only transverse double points. If $c$ is a double point of $\Pi_\C(L)$ then we write $\Pi_\C^{-1}(c)\cap L=\{c^+,c^-\}$, where $z(c^+)>z(c^-)$. (Since the Reeb field of the standard contact form $dz-y\,dx$ is $\pa_z$, there exists a 1--1 correspondence between double points of $\Pi_\C(L)$ and Reeb chords on $L$ and we will use these two notions interchangeably.)

Let $D_{m+k}$ be the unit disk $D$ in the complex plane $\C$ with punctures  $x_1,\dots,x_m,y_1,\dots,y_k$ on the boundary and let $J$ be an almost complex structure on $T^\ast M$ which is tamed by the standard symplectic form $\omega=dx\wedge dy$ on $T^\ast M$.
\begin{dfn}\label{d:phol}
A $J$-holomorphic disk with positive punctures $p_1,\dots,p_m$ and negative punctures $q_1,\dots,q_k$ and with boundary on $L$, is a map $u\colon D_{m+k}\to T^\ast M$ with the following properties.
\begin{itemize}
\item $\bar\pa_J u=du + J\circ du\circ i=0$ (where $i$ is the complex structure on the complex plane).
\item The restriction $u|\pa D_{m+k}$ has a continuous lift $\tilde u\colon\pa D_{m+k}\to L\subset J^1(M)$.
\item $\lim_{\zeta\to x_j}u(\zeta)=p_j$ and $\lim_{\zeta\to x_j\pm}{\tilde u}(\zeta)=p_j^{\pm}$, where $\lim_{\zeta\to x_j+}$ means that $\zeta$ approaches $x_j$ from the region in $\pa D_{m+k}$ in the positive direction as seen from $x_j$ and $\lim_{\zeta\to x_j-}$ means it approaches $x_j$ from the region in the negative direction.
\item $\lim_{\zeta\to y_j}u(\zeta)=q_j$ and $\lim_{\zeta\to y_j\pm}{\tilde u}(\zeta)=q_j^{\mp}$.
\end{itemize}
\end{dfn}
In this paper we will restrict attention to holomorphic disks with exactly one positive puncture. The moduli space of holomorphic disks with boundary on $L$ has certain compactness properties: a version of Gromov's compactness theorem holds, see \cite{BEHWZ, EES1, EES4}. Since the linearized $\bar\pa_J$-operator is elliptic and since the double points of $\Pi_\C(L)$ are transverse, the linearization of the equation which defines the moduli space of $J$-holomorphic disks is a Fredholm operator and its index determines the expected dimension of the moduli space. Details of the computation of the Fredholm index can be found in \cite{EES1, EES4}. Here we simply state the result. Pick for each Reeb chord $c$ of $L$ a capping path $\gamma_c$ in $L$ connecting $c^+$ to $c^-$. Define $\Gamma_c$ to be the path of Lagrangian subspaces $\Pi_\C(T_{\gamma_c} L)$ closed up by a positive rotation. Define
\[
|c|=\mu(\Gamma_c)-1,
\]
where $\mu$ is the Maslov index. Given any boundary condition of a holomorphic disk with positive puncture at $a$ and negative punctures at $b_1,\dots,b_k$ we can close it up to a loop by adding appropriately oriented capping paths at the Reeb chords. Let $A\in H_{1}(L;\Z)$ denote the homology class of the closed up loop and let $\M_A(a;b_1,\dots,b_k)$ denote the moduli space of holomorphic disks with positive puncture $a$, negative punctures at $b_1,\dots,b_k$, and boundary condition inducing the homology class $A$. Then
\[
\dim\bigl(\M_A(a;b_1,\dots,b_k)\bigr)=|a|-\sum_{j=1}^k|b_k|+\mu(\Gamma_A)-1,
\]
where $\Gamma_A$ is the path of Lagrangian subspaces along a loop in $L$ representing $A$.

\subsection{Legendrian contact homology}\label{s:LCH}
Let $L\subset J^1(M)$ be a Legendrian submanifold. The DGA $\A(L)$ of $L$ is the free algebra over the $\Z_2$-group ring of $H_1(L;\Z)$ generated by the Reeb chords of $L$,
\[
\A(L)=\Z_2[H_1(L;\Z)]\la c_1,\dots,c_m\ra.
\]
The differential $\pa\colon\A(L)\to\A(L)$ counts rigid holomorphic disks. It is defined to be linear over the algebra coefficients and to satisfy Leibniz rule and is thus determined by its action on generators. If $a$ is a generator then
\[
\pa a =\sum_{\dim(\M_A(a;b_1,\dots,b_k))=0}|\M_A(a;b_1,\dots,b_k)|\,A\,b_1\dots b_k,
\]
where $|\M|$ denotes the modulo $2$ number of points in the compact $0$-manifold $\M$. For this definition to make sense we require that $L$ is generic with respect to holomorphic disks so that all moduli spaces of holomorphic disks of dimension $\le 1$ are transversely cut out. It is shown in \cite{EES1, EES3} how to achieve such transversality by perturbing $L$. In these papers it is also shown that $\pa\colon \A(L)\to\A(L)$ is a differential (i.e., $\pa^2=0$) and that the {\em contact homology} $\krn(\pa)/\img(\pa)$ is invariant under Legendrian isotopies.

\begin{rmk}
In the case that $L$ is a spin manifold and $M$ is orientable one can lift the DGA over $\Z_2[H_1(L;\Z)]$  described above to a DGA over $\Z[H_1(L;\Z)]$. In this paper we will, as mentioned in Section \ref{S:intr}, concentrate on the $\Z_2$-case. This is mainly for simplicity, the results of the paper have straightforward generalizations to the more general setting of $\Z$-coefficients.
\end{rmk}

\subsection{Flow trees}\label{s:fltr}
We refer to \cite{E} for the detailed definitions of flow trees and only sketch the main points here. Let $L\subset J^{1}(M)$ be a Legendrian submanifold. Locally around points outside a codimension one subset of $M$, the image of $L$ under the {\em front projection} $\Pi_F\colon J^{1}(M)\to J^{0}(M)= M\times\R$, can be described as the graph of a finite number of functions. If $M$ is equipped with a Riemannian metric then these local functions define local gradients. We say that a curve $\gamma(t)$ in $M$ is a {\em flow line} of $L$ if it satisfies the differential equation
$$
\dot\gamma(t)=-\nabla(f_1-f_2)(\gamma(t)),
$$
where $f_1$ and $f_2$ are local functions of $L$. Any flow line has a natural {\em $1$--jet lift}, which consists of two curves in $L$ lying over $\gamma$. These curves are naturally oriented by the lifts of the vectors $-\nabla(f_1-f_2)$ and $-\nabla(f_2-f_1)$ to the sheets corresponding to $f_1$ and $f_2$, respectively. Projecting the $1$--jet lift to the cotangent bundle we get the {\em cotangent lift}.

\begin{dfn}\label{d:fltr}
A {\em flow tree of $L\subset J^1(M)$} is a continuous map $\phi\colon\Gamma\to M$, where $\Gamma$ is a source-tree which satisfies the following conditions.
\begin{itemize}
\item[{\rm (a)}] If $e$ is an edge of $\Gamma$ then $\phi\colon e\to M$ is an injective parametrization of a flow line of $L$.
\item[{\rm (b)}] Let $v$ be a $k$-valent vertex with cyclically ordered adjacent edges $e_1,\dots, e_k$. Let $\{\bar\phi_j^1,\bar\phi_j^2\}$ be the cotangent lift corresponding to $e_j$, $1\le j\le k$. We require that there exists a pairing of lift components such that for every $1\le j\le k$ (with $k+1=1$)
$$
\bar\phi_j^2(v)=\bar\phi_{j+1}^1(v)=\bar m\in\Pi_\C(L)\subset T^\ast M,
$$
and such that the flow orientation of $\bar\phi_j^2$ at $\bar m$ is directed toward $\bar m$ if and only if the flow orientation of $\bar\phi_{j+1}^1$ at $\bar m$ is directed away from $\bar m$.
\item[{\rm (c)}] The cotangent lifts of the edges of $\Gamma$ fit together to an oriented curve $\bar\phi$ in $\Pi_\C(L)$. We require that this curve is closed.
\end{itemize}
\end{dfn}
For simpler notation we will often denote flow trees simply by $\Gamma$, suppressing the parametrization map $\phi$ from the notation. We will also write $\bar\Gamma$ and $\tilde\Gamma$ for the cotangent and the $1$--jet lifts of $\Gamma$, respectively.

We next define punctures of a flow tree $\Gamma$. Let $v$ be a $k$-valent vertex of $\Gamma$ with cyclically ordered edges $e_1,\dots,e_k$. Consider two paired cotangent lifts $\bar\phi_j^2$ and $\bar\phi_{j+1}^1$ and the corresponding $1$--jet lifts $\tilde\phi_j^2$ and $\tilde\phi_{j+1}^1$ at $v$. If $\tilde\phi_j^2(v)\ne \tilde\phi_{j+1}^1(v)$ then both must equal Reeb chord endpoints. If this is the case then we say that $v$ contains a {\em puncture} after $e_j$. It is shown in \cite{E} that any flow tree with a vertex $v$ which contains more than one puncture is a union of flow trees such that every vertex of each one of them contains at most one puncture. Thus we may restrict attention to flow trees with at most one puncture at each vertex and we call such a vertex a {\em puncture} of the tree.

Let $p$ be a puncture of a flow tree. Let $\tilde\phi^1$ and $\tilde\phi^2$ be the $1$--jet lifts which map to the Reeb chord at $p$, with notation chosen so that $\tilde\phi^1$ is oriented toward $\tilde\phi^1(p)$ and $\tilde\phi^2$ oriented away from $\tilde\phi_2(p)$. Then we say that $p$ is a {\em positive} puncture if
$$
z(\tilde\phi^1(p))<z(\tilde\phi^2(p)),
$$
and we say it is {\em negative} is the opposite inequality holds. (Recall, $z$ is the coordinate in the $\R$-direction of $J^1(M)=T^\ast M\times\R$.) Using the symplectic area of a flow tree it is not hard to show that every flow tree $\Gamma$ in $M$ has at least one positive puncture. If $q$ is a puncture of a flow tree then $q$ corresponds to some Reeb chord of $L$ which in turn corresponds to a critical point of some local function difference of $L$. If $L$ is generic then the Hessian at this critical point is non-degenerate. Let $I(q)$ denote the index of the critical point of the positive function difference corresponding to the Reeb chord.

As for holomorphic disks there is a simple dimension formula for flow trees which computes the expected dimension of the space of flow trees with $1$--jet lift homotopic to that of the given tree. For trees with exactly one positive puncture this formula reads
$$
\dim(\Gamma)=(I(p)-1)
-\sum_{q\in Q(\Gamma)}(I(q)-1)
+\sum_{r\in R(\Gamma)}\mu(r)-1,
$$
where $p$ is the positive puncture of $\Gamma$, $Q(\Gamma)$ its set of negative punctures, $R(\Gamma)$ its set of vertices which are not punctures, and $\mu(r)$ is the Maslov content of $r$, see \cite[Section 3.1.1]{E}. Viewing the $1$--jet lift as a boundary condition for a holomorphic disk the dimension formula just stated agrees with the dimension formula for holomorphic disks. In fact for rigid disks and trees more is true as the following theorem shows. (See \cite[Theorem 1.1]{E} for a more general version.)
\begin{thm}\label{t:disktree}
If $L\subset J^1(M)$ is a Legendrian submanifold of dimension $\le 2$ then there is a generic complex structure $J$ on $T^\ast M$ such that there is a 1--1 correspondence between rigid $J$-holomorphic disks with boundary on $L$ and rigid flow trees determined by $L$.
\end{thm}
This theorem has the consequence that one may replace holomorphic disks in the definition of Legendrian contact homology differential with flow trees. We will often do so below.

We will sometimes talk about convergence of flow trees. When doing so we employ the following topology on the space of flow trees. Let $\Gamma$ be a flow tree with one positive puncture $p_0$. Orient $\Gamma$ by declaring any edge starting at $p_0$ to be oriented away from $p_0$ and by requiring that at each vertex different from $p_0$ there is exactly one edge oriented toward it. We will associate a planar domain $\Delta_\Gamma$ to $\Gamma$ as follows. For each $1$-valent vertex different from $p_0$ and for each puncture different from $p_0$ pick a half infinite strip of width $1$. For each vertex where such edges begin, pick a strip of width which equals the sum of the widths of the edges going out from it and of length agreeing with the natural flow length of the edge coming in to that vertex. Gluing the half strips to the strip, we obtain a half strip domain of finite width with slits toward $+\infty$. (Here we consider the half strip of a negative puncture at a vertex of valence $\ge 2$ as outgoing from that vertex.) Continuing inductively in this way we get a strip domain with slits toward $+\infty$ and where $-\infty$ corresponds to $p_0$. This domain $\Delta_\Gamma$ is determined up to over all translation of the slits. Thus the conformal structure on $\Delta_\Gamma$ thought of as a disk with boundary punctures is uniquely determined by $\Gamma$. 

Note that the cotangent lift of $\Gamma$ is naturally defined as a map from the boundary $\pa\Delta_\Gamma$ to $\Pi_\C(L)$. Consider the bundle of mapping spaces of maps from the boundary of a punctured disk into $\Pi_\C(L)$ over the space of conformal structures on the disk with $m$ punctures, endowed with the $C^0$-topology. We topologize the space of flow trees by giving it the subspace topology with respect to this topology on the bundle of mapping spaces.
\section{Generalizations of Legendrian contact homology}\label{S:gench}
In this section, we generalize Legendrian contact homology slightly. The first generalization concerns a special kind of non-closed Legendrian submanifolds which can be concatenated in a natural way, see Subsection \ref{s:ends}. The second generalization, see Subsections \ref{s:genfltr} and \ref{s:pert},  concerns stabilization of Legendrian submanifolds. It is more elaborate than the first and is inspired by Morse-Bott techniques for contact homology, see \cite{B}, as well as so called abstract perturbations, see \cite{FOOO, HWZI, HWZII}.

\subsection{Holomorphic disks and flow trees of $1$-dimensional Legendrian submanifolds}\label{s:1Ddisktree}
Let $L\subset J^1(\R)$ be a Legendrian submanifold. Using the Riemann mapping theorem, moduli-spaces of holomorphic disks with boundary on $L$ can be understood geometrically as follows. The Reeb chords of the diagram correspond to double point of the Lagrangian projection $\Pi_\C\colon J^1(\R)\to T^\ast\R\approx \C$ and a holomorphic disk to an ordinary holomorphic disk with boundary on the knot diagram and such that at the positive puncture the incoming strand of the disk is upper and at a negative puncture the incoming strand is lower.

\begin{lma}\label{l:homeo}
Let ${\mathcal F}$ denote the space of flow trees with one positive puncture determined by $L$, let $\M$ denote the moduli space of all holomorphic disks with one positive puncture determined by $L$ and let $\M'\subset \M$ denote the subspace of all disks without interior branch points.
Then the spaces $\mathcal{F}$ and $\M'$ are homeomorphic and the subset $\M'\subset \M$ is closed and has a natural compactification consisting of broken disks from $\M'$.
\end{lma}

\begin{pf}
Let $\Gamma$ be a tree. Joining corresponding points on the cotangent lift of $\Gamma$ by straight lines in the fibers of $T^{\ast}\R$, we obtain a map of a disk with branch points on the boundary. By the Riemann mapping theorem it admits a holomorphic parametrization. Conversely, any holomorphic disk without interior branch points gives flow lines joined at branch points. That is, a flow tree with higher valence vertices corresponding to branch points. Thus there is a natural bijection between ${\mathcal F}$ and $\M$ as sets. We then need only show that the map $\M\to{\mathcal F}$ is continuous and open. This however is straightforward: a holomorphic disk is determined by the location of its branch points, varying these give an open subset in the space of trees. Finally, it is clear that the subset of the compactified moduli space of holomorphic disks which consists of (broken) disks without interior branch points is closed. Being a closed subset of a compact space it is compact as well.
\end{pf}

\begin{cor}\label{c:compf}
The space of flow trees is a manifold with boundary with corners. In particular, the boundary of a component of the space of trees consists of broken trees.
\end{cor}

\begin{pf}
This holds for the compactified moduli space of holomorphic disks. The corollary thus follows from Lemma \ref{l:homeo}.
\end{pf}

\begin{rmk}
In order to see the relation between $\M$ and $\M'$, note that any disk with a interior branch points can be connected to a disk of the same dimension with all its branch points on the boundary by pushing the branch points. For example, a second order branch point which is pushed to the boundary becomes a third order boundary branch point which splits into two ordinary boundary branch points. A local model can be obtained as follows. Consider the map $f(z)=z^2$ with domain $\Omega$ bounded by the curve $\{2xy=-\epsilon, x>0\}$, $\epsilon>0$, and containing the coordinate axes. Note that $f$ maps $\pa\Omega$ to the line $y=-\epsilon$ and that it has a branch point at $0$. To see what happens as the branch point approaches the boundary we let $\epsilon\to 0$. In the limit $f$ becomes a map $f_0$ from the complement of the open fourth quadrant of the complex plane. Consider the map $z(w)=w^{3/2}$, taking the upper half plane to the domain of $f_0$. Then the map $f_0(z(w))$ looks like $w\mapsto w^3$ near the origin and the $2$--parameter family continues as $w\mapsto \delta_1+w(w^2-\delta_2)$ for $\delta_1,\delta_2\in\R$, $\delta_2>0$.
\end{rmk}

\subsection{Legendrian submanifolds with standard ends}\label{s:ends}
Let $L\subset J^1(\R^n)$ be a Legendrian submanifold parametrized by
$\gamma(p)=(x(p),y(p),z(p))\in J^1(\R^n)$. For $\kappa>0$, consider the Legendrian submanifold $\tilde L[\kappa,\tau]\subset J^1(\R^{n+1})$ parametrized by
$$
\tilde\gamma[\kappa,\tau](p,t)=\Bigl(t,x(p),2\kappa(t-\tau)z(p),
\kappa\bigl(1+(t-\tau)^2\bigr)\cdot y(p),
\kappa\bigl(1+(t-\tau)^2\bigr)\cdot z(p)\Bigr).
$$
Note that the Reeb chords of $\tilde L[\kappa,\tau]$ correspond to the Reeb chords of $L$, lying in the slice $x_0=\tau$.

Let $M$ be an $(n+1)$-manifold with boundary $\pa M= L_-\sqcup L_+$. Let $f\colon M\to J^1(\R^{n+1})$ be a Legendrian embedding such that $|x_0\circ f|\le 1$, such that $|x_0\circ f|<1$ in the interior of $M$,  and such that $f$ agrees with $\gamma_+[\kappa,+1]\colon L_+\times(-\epsilon,0]$ in some collar neighborhood $L_+\times(-\epsilon,0]$ of $L_+$, and such that it agrees with $\gamma_-[\kappa,-1]\colon L_-\times[0,\epsilon)$ in some collar neighborhood $L_-\times[0,\epsilon)$ of $L_+$, where $\gamma_\pm$ are some Legendrian embeddings of $L_\pm$. We say that $f(M)$ is a {\em Legendrian embedding with standard ends}.

If $f\colon M\to J^1(\R^{n+1})$ is a Legendrian embedding with standard ends, then we can construct a Legendrian embedding $F$ of a non-compact manifold $\hat M$ obtained by adding $\gamma_+(L_+\times[0,\infty))$ and $\gamma_-(L_-\times(-\infty,0])$ to $F$. Note that the Reeb chords of $F$ are exactly those of $f$. Moreover we have the following.

\begin{lma}\label{l:ends+disk}
Any holomorphic disk of finite area with boundary on $F(\hat M)$ lies in the region $|x_0|\le 1$ and the space of holomorphic disks is Gromov compact. Moreover any holomorphic disk with its positive puncture at $x_0=\pm 1$ lies entirely in this slice.
\end{lma}

\begin{pf}
As in \cite{EES3}, projection to the $(x_0+iy_0)$-line in combination with the maximum principle establishes the first and third statements. The second then follows from the standard proof of Gromov compactness, see e.g. \cite{EES1}.
\end{pf}

We have the corresponding statement for flow trees.
\begin{lma}\label{l:ends+trees}
Any flow tree of $F(\hat M)$ of finite symplectic area is contained in the region $|x_0|\le 1$ and any flow tree with its positive puncture at $x_0=\pm 1$ stays entirely in this slice. Moreover, there is a $1-1$ correspondence between on the one hand rigid holomorphic disks with boundary on $F(\hat M)$, and on the other rigid flow trees determined by $F(\hat M)$ and rigid flow trees of $L_\pm$ in slices.
\end{lma}

\begin{pf}
The statements about flow trees follows from the fact that all gradient differences in the slices $x_0=\pm 1$ have trivial $x_0$-component. The second statement follows from a slight modification of the proof of \cite[Theorem 1.1]{E}, which consists of two parts: convergence of holomorphic disks to flow trees and construction of holomorphic disks near rigid flow trees. The convergence part of the proof holds without change in this more general setting. The construction part can be subdivided into two parts: for trees not in the slices $x_0=\pm 1$ the proof from \cite{E} applies. For trees in slices $x_0=\pm 1$ we apply the proof from \cite{E} to $L_\pm$ thought of as lying in the slice.
\end{pf}

Lemma \ref{l:ends+disk} implies that the contact homology of $F(\hat M)$ is well defined. (Here we restrict attention to Legendrian isotopies in the class of Legendrian submanifolds with standard ends). Lemma \ref{l:ends+trees} implies that we can compute it using flow trees instead of holomorphic disks. We then define the contact homology of $f\colon M\to J^1(\R^{n+1})$ to equal the contact homology of $F\colon \hat M\to J^1(\R^{n+1})$.

\begin{rmk}\label{r:concatenation}
If $f_j\colon M_j\to J^{1}(\R^{n+1})$ is a Legendrian submanifold with standard ends and if one of the ends of $(f_0,M_0)$ agrees with an end of $(f_1,M_1)$ then the two Legendrian embeddings can be joined to a Legendrian embedding $f_{01}$ of the manifold $M_{01}$ obtained by joining $M_0$ and $M_1$ along their common boundary. The contact homology differential of the join is determined in a straightforward way by the contact homology differentials of its pieces.
\end{rmk}

\subsection{Generalized flow trees}\label{s:genfltr}
In this subsection we introduce the notion of generalized flow trees for the product of a given Legendrian submanifold $L\subset J^{1}(\R)$ and a manifold with boundary equipped with a Morse function. In the present paper we will apply this construction only in the case when the auxiliary manifold factor is an interval or a $2$-disk.

Let $N$ be a manifold with boundary $\pa N$. We will use Morse functions $\beta\colon N\to\R$ of the following form. The restriction $\beta_{\pa}=\beta|\pa N$ is a Morse function on $\pa N$ and there is a collar neighborhood $\pa N\times [0,\epsilon)$ of $\pa N$ in $N$ where $\beta(x,t)= \beta_\pa(x) + k t^{2}$, where $(x,t)\in\pa N\times[0,\epsilon)$ and where $k>0$ is a constant. We call a Morse function of this type {\em boundary adjusted}.

Consider the Legendrian submanifold $L\times N\subset J^{1}(\R\times N)$. The Reeb chords of this Legendrian submanifold come in $N$-families, one for each Reeb chord of $L$. We denote the manifold of Reeb chords corresponding to the Reeb chord $c$ by $N_c$. Choose a boundary adjusted Morse function $\beta_c\colon N_c\to\R$ for each Reeb chord manifold $N_c$. Endow $N$ with a Riemannian metric $g$ which has the following form in the collar neighborhood of the boundary
\begin{equation}\label{e:bdrymetric}
g(x,t)=g_\pa(x) + dt^{2},
\end{equation}
for $(x,t)\in\pa N\times[0,\epsilon)$, where $g_\pa$ is a Riemannian metric on $\pa N$.

A {\em flow line} $\gamma$ in a Reeb chord manifold $N_c$ is an oriented segment which can be parametrized in an orientation preserving manner by a solution to the gradient equation
$$
\dot q=-\nabla\beta_c(q),\quad q\in N,
$$
where the gradient is defined using the Riemannian metric $g$ of \eqref{e:bdrymetric}. Note that a gradient line which starts in $\pa N_c$ stays in $\pa N_c$ and that a gradient line which starts in $N_c-\pa N_c$ can hit $\pa N_c$ only at a critical point in the boundary.

In order to define the notion of a generalized tree we first introduce some preliminary concepts.
If $\Gamma$ is a flow tree of $L$ and if $q_0\in N$ then we let $\Gamma_{q_0}$ denote this flow tree considered as a flow tree of $L\times N$ and lying in the slice $\{(q,x)\in N\times\R\colon q=q_0\}$. We call $\Gamma_{q_0}$ a {\em slice tree}. A {\em level} is a finite collection $\Theta=\{\Gamma_{q_1},\dots,\Gamma_{q_r}\}$ of (unbroken) slice trees. A {\em connector} is a finite collection of flow lines $\theta=\{\gamma_1,\dots,\gamma_k\}$ in Reeb chord manifolds, where we allow also flow lines of length $0$.

A {\em generalized flow tree with $k$ levels} is an ordered collection $(\Theta^1,\dots,\Theta^k)$ of levels together with an ordered collection $(\theta^0,\dots,\theta^k)$ of connectors which have the following properties.
\begin{itemize}
\item The connector $\theta^0$ consists of exactly one flow line in $N_a$ emanating at a critical point $p^0\in N_a$ of $\beta_a$ for some Reeb chord $a$ and ending at $q^0\in N_a$. If $k=0$ then $q^0$ is a critical point of $\beta_\alpha$ as well, if $k>0$ then $q_0$ is not a critical point. Let $p^1=q^0$.
\item The level $\Theta^1$ consists of exactly one slice tree $\Gamma_{p^1}$ with positive puncture at $a$.
\item Let $j>0$ and let $(\Gamma^{j,1})_{p^{j}_1},\dots,(\Gamma^{j,r})_{p^{j}_r}$ denote the slice trees in the $j^{\rm th}$ level $\Theta^{j}$ and let $c^{j,s}_{1},\dots,c^{j,s}_{m_s}$ denote the Reeb chords at the negative punctures of $\Gamma^{j,s}$, $s=1,\dots,r$. Then $\theta^j$ consists of flow lines emanating from those $p^{j}_s\in N_{c^{j,s}_t}$, $s=1,\dots, r;t=1,\dots,m_s$ which are not critical points of $\beta_{c^{j,s}_t}$. Let $q^{j}_1,\dots,q^{j}_l$ be the endpoints of the flow lines $\gamma^{j}_1,\dots,\gamma^{j}_l$ in $\theta^{j}$ which are not critical points. Let $p^{j+1}_s=q^{j}_s$, $s=1,\dots,l$.
\item For $0<j\le k$, $\Theta^{j+1}$ consists of slice trees $(\Gamma^{j+1,1})_{p^{j+1}_1},\dots,(\Gamma^{j+1,r})_{p^{j+1}_l}$ such that the positive puncture of $\Gamma^{j+1,s}$ is at the Reeb chord $a^{j+1}_s$ where $\gamma^{j}_s\subset N_{a^{j+1}_s}$, $s=1,\dots,l$.
\item Let $q\in N_c$ be a point where some flow line $\gamma$ in a connector $\theta^{s}$, $s=0,1,\dots,k$ ends. Then if there is no level tree $\Gamma_q$ in $\Theta^{s+1}$ with positive puncture matching the Reeb chord $c$ then $q$ is a critical point of $\beta_c$.
\end{itemize}

If $G$ is a generalized flow tree as just described then we say that $G$ has positive puncture at the critical point of $a^{j}\in N_a$ of $\beta_a$ and negative punctures at all negative punctures of slice trees in $\Theta^j$ where no flow line in $\theta^j$ begins (these are critical points $b_s^{t}\in N_{b_s}$) and at the critical points $b_s^{t}\in N_{b_s}$ where some flow lines in some $\theta^l$, $1\le l\le k$, ends. Any generalized flow tree has natural $1$--jet and cotangent lifts. Adding suitably oriented capping paths at the punctures to the projection of the $1$--jet lift of $G$ to $J^{1}(\R)$ we get a homology class $A\in H_1(L;\Z)$. We write $\Co_A(a^{j};b_1^{j_1},\dots b_m^{j_m})$ for the space of such generalized trees.

As for flow trees we associate a planar domain and a map of its boundary into $T^\ast\R^2$ to any generalized flow tree. Using the map from the boundary of the planar domain we may endow the set of generalized flow trees with a topology. It is a straightforward consequence of the compactness properties of flow trees, see Lemma \ref{c:compf}, that also the space of generalized flow trees has a natural compactification consisting of broken generalized flow trees. We denote the space of generalized flow trees by $\Co$ and its compactification ${\overline\Co}$.

\subsection{Perturbed generalized flow trees}\label{s:pert}
Let $L\subset J^{1}(\R)$ be a Legendrian submanifold and let $N$ be a manifold with boundary. As in Subsection \ref{s:genfltr} we consider the Legendrian submanifold $L\times N\subset J^{1}(\R\times N)$ and we equip $N$ with a Riemannian metric and each Reeb chord manifold $N_c$, $c\in\Rch(L)$, with a Morse function $\beta_c$ satisfying conditions as stated there. Let ${\overline \M}$ denote the moduli space of flow trees on $L$ and fix a function $v\colon{\overline \M}\times N\to C^{k}(S^{1},TN)$ with the following properties.
\begin{itemize}
\item $v(\Gamma,n)\in T_n N$ and $v(\Gamma,n)$ is tangent to $\pa N$ for $n\in\pa N$.
\item By scaling lengths, we think of the source $S^{1}$ of $v(\Gamma,n)$ as the cotangent lift $\bar \Gamma$ of $\Gamma$ (with resolved self intersections). We require $v$ to be constant in neighborhoods of the punctures of $\bar\Gamma$ and equal to zero near the positive puncture. Consider a broken tree $\Gamma$ containing a broken tree $\Gamma''$. Let $\Gamma'$ denote the broken tree obtained by removing $\Gamma''$ from $\Gamma$ and assume that $\Gamma''$ is attached to $\Gamma'$ at a point $y$ in its cotangent lift. Then we require that for $x\in\Gamma''$ the following {\em join equation} holds
\begin{equation}\label{e:join}
v(\Gamma,n)(x)=v(\Gamma',n)(y)+v(\Gamma'',n)(x).
\end{equation}
\end{itemize}
We call a function $v\colon {\overline M}\times N\to C^{k}(S^{1},TN)$ with these properties a {\em perturbation function}.

The definition of a perturbed flow tree is analogous to that of a generalized flow tree. Fix a perturbation function $v$. If $\Gamma$ is a flow tree of $L$ and if $q\in N$ then recall that $\Gamma_q$ was used to denote the slice tree corresponding to $\Gamma$. Think $\Gamma_q$ as a map $\Gamma_q=(\Gamma^{\R},q)\colon S^{1}\to \R\times N$ where $S^1$ is the cotangent lift of $\Gamma$ and the map $(\Gamma^{\R}(t),q)$ equals the natural map into $T^{\ast}\R\times T^{\ast} N$ followed by projection to $\R\times N$. We let $\tilde\Gamma_{q}$ denote the following map $S^{1}\to\R\times N$
\[
\tilde\Gamma_q(t)=\left(\Gamma^{\R}(t),\exp_q\bigl(v(\Gamma,q)(t)\bigr)\right),
\]
where $\exp$ denotes the exponential map in a Riemannian metric of the form given in \eqref{e:bdrymetric}.
We call $\tilde\Gamma_{q}$ a {\em perturbed slice tree} and if $x$ is a negative puncture of $\Gamma$ then we write $e(q,x)=\exp_{q}(v(\Gamma,q)(x))$. A {\em perturbed level} is a finite collection $\tilde\Theta=\{\tilde\Gamma_{q_1},\dots,\tilde\Gamma_{q_r}\}$ of (unbroken) perturbed slice trees. A {\em connector} is a finite collection of flow lines $\theta=\{\gamma_1,\dots,\gamma_k\}$ in Reeb chord manifolds, where we allow also flow lines of length $0$.

A {\em perturbed flow tree with $k$ levels} is an ordered collection $(\tilde\Theta^1,\dots,\tilde\Theta^k)$ of levels together with an ordered collection $(\theta^0,\dots,\theta^k)$ of connectors which have the following properties.
\begin{itemize}
\item The connector $\theta^0$ consists of exactly one flow line in $N_a$ emanating at a critical point $p^0\in N_a$ of $\beta_a$ for some $a\in\Rch(L)$ and ending at $q^0\in N_a$. If $k=0$ then $q^0$ is a critical point of $\beta_\alpha$ as well, if $k>0$ then $q_0$ is not a critical point. Let $p^1=q^0$.
\item The level $\tilde\Theta^1$ consists of exactly one slice tree $\tilde\Gamma_{p^1}$ with positive puncture at $a$.
\item Let $j>0$ and let $(\tilde\Gamma^{j,1})_{p^{j}_1},\dots,(\tilde\Gamma^{j,r})_{p^{j}_r}$ denote the slice trees in the $j^{\rm th}$ level $\tilde\Theta^{j}$ and let $c^{j,s}_{1},\dots,c^{j,s}_{m_s}$ denote the Reeb chords at the negative punctures of $\Gamma^{j,s}$, $s=1,\dots,r$. Then $\theta^j$ consists of flow lines emanating from those $e(p^{j}_s,x)\in N_{c^{j,s}_t}$, $s=1,\dots, r;t=1,\dots,m_s$, where $x$ denotes the negative puncture corresponding to the Reeb chord $c^{j,s}_t$ which are not critical points of $\beta_{c^{j,s}_t}$. Let $q^{j}_1,\dots,q^{j}_l$ be the endpoints of the flow lines $\gamma^{j}_1,\dots,\gamma^{j}_l$ in $\theta^{j}$ which are not critical points. Let $p^{j+1}_s=q^{j}_s$, $s=1,\dots,l$.
\item For $0<j\le k$, $\tilde\Theta^{j+1}$ consists of perturbed slice trees $(\tilde\Gamma^{j+1,1})_{p^{j+1}_1},\dots,(\tilde\Gamma^{j+1,r})_{p^{j+1}_l}$ such that the positive puncture of $\Gamma^{j+1,s}$ is at the Reeb chord $a^{j+1}_s$ where $\gamma^{j}_s\subset N_{a^{j+1}_s}$, $s=1,\dots,l$.
\item Let $q\in N_c$ be a point where some flow line $\gamma$ in a connector $\theta^{s}$, $s=0,1,\dots,k$ ends. If there is no perturbed level tree $\tilde\Gamma_q$ in $\tilde\Theta^{s+1}$ with positive puncture matching the Reeb chord $c$ then $q$ is a critical point of $\beta_c$.
\end{itemize}

In analogy with generalized flow trees, we use the notion $\Co_{{\rm p};A}(a^{j};b_1^{j_1},\dots b_m^{j_m})$ to denote the space of perturbed trees with specified punctures and boundary data as well as the notions $\Co_{\rm p}$ for the space of all perturbed flow trees and ${\overline \Co}_{\rm p}$ for its natural compactification consisting of several level trees.

Given a perturbation function and Morse functions on all Reeb chord manifolds we associate a contact homology algebra with a differential to the Legendrian submanifold $L\times N$. The algebra $\A(L\times N)$ is the free $\Z_2[H_1(L;\Z)]$-algebra generated by critical points of the Morse functions $\beta_c$, for Reeb chords $c$ and the differential $\pa_{\rm p}\colon\A\to\A$ counts rigid perturbed flow trees. More precisely it satisfies Leibniz rule and is defined as follows on generators,
\[
\pa_{\rm p} a = \sum_{\dim(\Co_{{\rm p};A}(a;b_1,\dots,b_k))=0}|\Co_{{\rm p};A}(a,b_1,\dots,b_k)|A\,b_1\dots b_k.
\]

\begin{lma}\label{l:pertalg}
For generic sufficiently small perturbations and Morse functions $\beta_c$, the moduli space of perturbed flow trees of formal dimension $\le 1$ are transversely cut out. Consequently, $\pa_{\rm p}$ is a differential i.e. $\pa^{2}_{\rm p}=0$.
\end{lma}

\begin{pf}
The transversality properties follows from standard applications of the finite dimensional jet-transversality theorem. Furthermore it is clear that, for small enough perturbation function, any sequence of perturbed rigid trees converges to a possibly broken perturbed rigid tree. The statement on the differential then follows from the usual gluing argument, which is technically very easy in this case since we need only glue Morse flow lines. (The perturbation function has to be sufficiently small so that for any $(\Gamma,n)\in\M\times N$, $v(\Gamma,n)$ is smaller than the injectivity radius of the Riemannian metric.)
\end{pf}

We next concentrate on the case important for our applications. Let $N=I=[-1,1]$ and consider $L\times I$ equipped with two different perturbation functions $v^{0}$ and $v^{1}$ and with two different sets of Morse functions $\{\beta_c^{0}\}$ and $\{\beta_c^{1}\}$ where $c$ ranges over all Reeb chords. These data determines two different differentials $\pa_0$ and $\pa_1$ on $\A(L\times I)$. We show next that the resulting DGAs are tame isomorphic. We restrict the Morse functions so that they have minima at $\pm 1$ and exactly one interior maximum. We then write the algebra generators as $\{\hat c, c[-1], c[+1]\}$, where $\hat c$ denotes the maximum in $I_c$,  $c[\pm 1]$ the minima at $\pm 1\in N_c$, and where $c$ ranges over all Reeb chords.

In order to see how the differentials $\pa_0$ and $\pa_1$ are related, we chose a generic path $(v^{s},\{\beta^{s}_c\})$, $0\le s\le 1$ of perturbations and Morse functions connecting the two sets of given data, where we restrict the Morse functions to be of the form described above. With such a path chosen we get, for each $s$ a moduli spaces $\Co_{\rm p}^{s}$ of perturbed trees of $L\times I$.

\begin{lma}\label{l:(-1)tree}
For a generic path there exists no perturbed flow trees in $\Co_{\rm p}^{s}$ of formal dimension $\le -2$ and there exists finitely many instances where there are perturbed flow trees of formal dimension $-1$ and at such an instance there is exactly one perturbed flow tree of dimension $-1$ which is transversely cut out in the sense of $1$-parameter families.
\end{lma}

\begin{pf}
This follows since the slice trees are transversely cut out and hence the only degeneration possible is when some incidence equation (involving the perturbation function and stable/unstable manifolds of the Morse functions) has non-transverse solutions.
\end{pf}

Thus, for $s$ which is not a $(-1)$-tree instance, we get an induced differential $\pa_s$ on $\A(L\times I)$ by counting rigid perturbed trees.

\begin{rmk}
We note that the appearance of a $(-1)$-tree as in Lemma \ref{l:(-1)tree} implies that many of the perturbed flow trees of formal dimension $0$ are in fact appearing in higher dimensional families. These are obtained from rigid flow trees with some negative punctures at the positive puncture of the $(-1)$-disk by gluing. One way to relate the differentials is to study the details of how such families split as the deformation variable changes. Below we will however use a less explicit but technically simpler method.
\end{rmk}

To show that the DGAs $(\A(L\times I),\pa_0)$ and $(\A(L\times I),\pa_1)$ are tame isomorphic, we use a stabilization argument which is very close to the argument given in \cite[Section 4.3]{EES3}. Consider first the $1$-parameter family of moduli spaces $\Co_{{\rm p}}^{s}$ and let $\Co_{{\rm p}}^{s}(0)$ denote the part of $\Co_{{\rm p}}^{s}$ which consist of all perturbed trees of formal dimension $0$. If $s$ varies in a interval $[s_0,s_1]$ which does not contain any $(-1)$-tree instance then the $1$-manifold
$$
\cup_{s\in[s_0,s_1]}\Co_{{\rm p}}^{s}(0)
$$
gives a cobordism between $\Co_{{\rm p}}^{s_0}(0)$ and $\Co_{{\rm p}}^{s_1}(0)$ and consequently $\pa_{s_0}=\pa_{s_1}$.

\begin{figure}
\centering
\includegraphics[width=.4\linewidth]{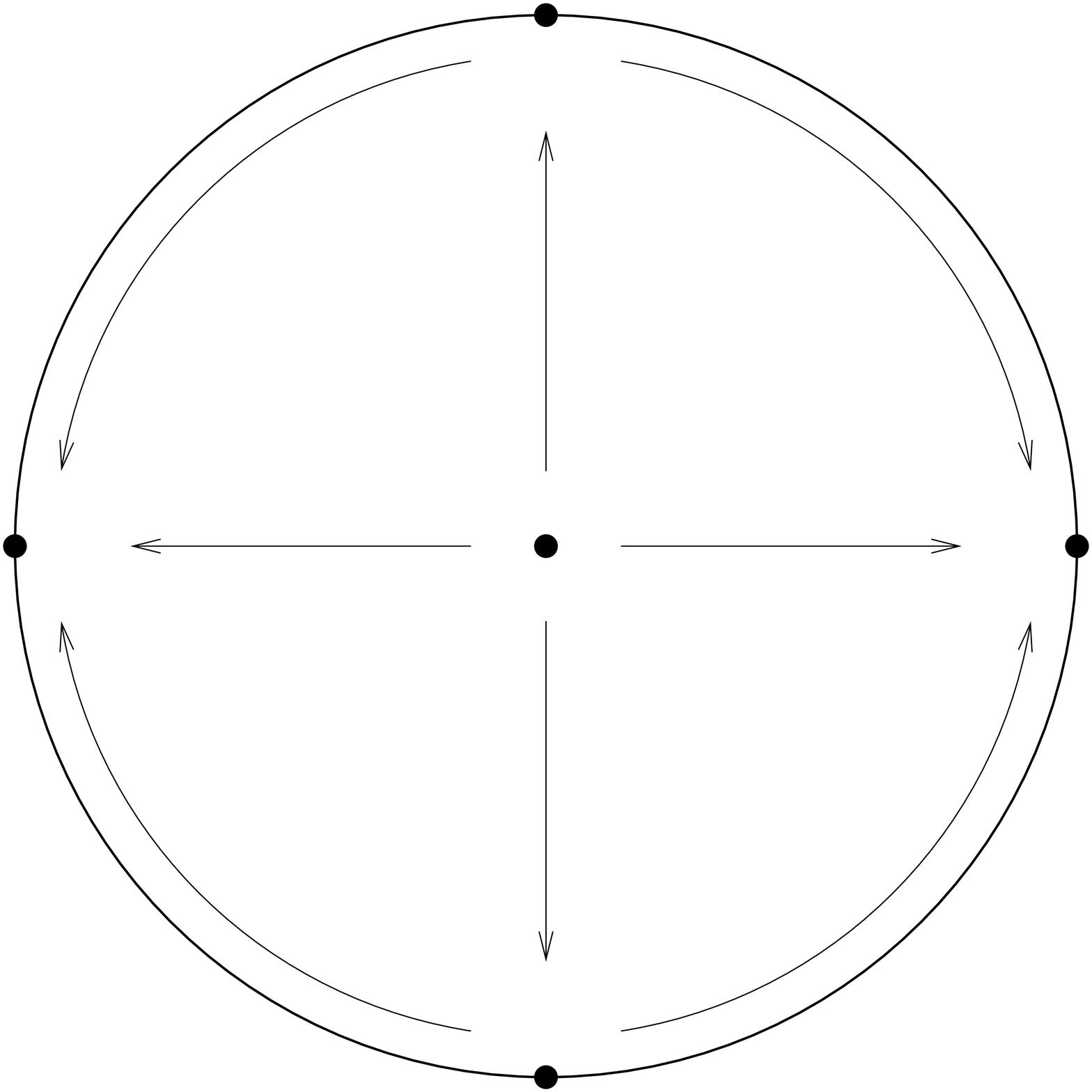}
\caption{The Morse flow on the disk. Flow lines connecting $\tilde c$ to $\hat c[\sigma]$ and $c[\sigma]$, as well as $\hat c[\sigma]$ to $c[\sigma]$ are indicated.}
\label{fig:disk}
\end{figure}

To show that $\pa_{0}$ and $\pa_{1}$ are tame isomorphic it is thus sufficient to show that the DGAs on both sides of a $(-1)$-disk instance are tame isomorphic. We consider the product $L\times D$ where $D$ is a $2$-disk. Think of the boundary of this disk as consisting of $I_0$ and $I_1$. Let $\beta_{c}^{0}$ and $\beta_{c}^{1}$ denote the Morse functions on the Reeb chord manifolds on the two sides of the $(-1)$-disk moment. Choose the Morse functions $\beta_c\colon D\to\R$ as extensions of these functions with exactly one maximum and with Morse flows as shown in Figure \ref{fig:disk}.

Let furthermore the extension of the perturbation to the disk be that of the generic family, let  $\tilde\A=\A(L\times D)$, and let $\Delta$ denote the corresponding differential. Write $\tilde c$ for the maximum Reeb chords, $\hat c[0]$ and $\hat c[1]$ for the saddle Reeb chords and $c[0]$ and $c[1]$ for the minimum Reeb chords. Let $\Omega\colon \tilde\A\to\A$ be the algebra morphism which takes $\tilde c$ to $0$, $\hat c[\sigma]$ to $\hat c$ and $c[\sigma]$ to $c$.

\begin{lma}\label{l:diffconst}
The differential $\Delta$ satisfies
\begin{align}\label{e:tame}
&\Delta c[\sigma]=(\pa_1 c)[\sigma]=(\pa_0 c)[\sigma],\quad \sigma=0,1,\\\notag
&\Delta \hat c[\sigma]=(\pa_\sigma \hat c)[\sigma],\quad \sigma=0,1,\\\notag
&\Delta \tilde c =\hat c[0]+\hat c[1] + \epsilon + \Ordo(1),
\end{align}
where $\Omega\epsilon=m\alpha$, where $\alpha$ is the word of negative punctures of the $(-1)$-disk multiplied by the homology class of its boundary condition, where $m\in\Z_2$, and where $\Ordo(1)$ denotes words which contain at least one $\tilde c$-variable.
\end{lma}

\begin{pf}
The first two equations follow from the nature of the metric and the Morse function near the boundary. The second property follows from a limiting argument where we let the disk limit to the middle family containing the $(-1)$-disk. Any rigid disk must limit to a disk of formal dimension $0$. All such disks which are not the $(-1)$-disk have at least one negative puncture at a $\tilde c$ variable.
\end{pf}

\begin{cor}\label{c:tameiso}
There exists a tame isomorphism between the algebras $(\A,\pa_0)$ and $(\A,\pa_1)$.
\end{cor}

\begin{pf}
This follows from the fact that $\Delta^{2}=0$ in combination with \eqref{e:tame}. The details of the proof are found in \cite[Lemma 4.21]{EES2}.
\end{pf}

\section{Geometric perturbations of the trace of the constant isotopy}\label{S:geomprt}
In this section we study the contact homology of Legendrian submanifolds which are traces of constant isotopies. We show that the contact homology differential can be expressed in terms of perturbed flow trees in this case.

\subsection{Structure of the perturbation}\label{s:strprt}
Let $L\subset J^1(\R)$ be a Legendrian submanifold with Reeb chords $c_1,\dots,c_m$ and let $L\times I\subset J^1(\R^2)$ denote the trace of the constant isotopy starting and ending at $L$. For convenient notation we take $I=[-1,1]$ and we think of $L\times I\subset L\times\R$. Note that the Reeb chords of $L\times\R$ form $1$-dimensional manifolds $\approx [-1,1]$, one for each Reeb chord of $L$. We will denote these manifolds $M_c\subset T^\ast \R^2$ where $c\in\{c_1,\dots, c_m\}$ is the corresponding Reeb chord of $L$ and we will sometimes think of these manifolds as submanifolds of $J^0(\R^2)$. Furthermore we will write $M_c^+$ and $M_c^-$ for the corresponding submanifolds of $L\times I$ of Reeb chord endpoints. In order to define the Legendrian contact homology of $L\times I$ we use any small perturbation which turns $L\times I$ into a Legendrian submanifold with standard ends and employ the definition from Subsection \ref{s:ends}.

To get a more detailed understanding of the contact homology of $L\times I$ we will choose very specific perturbations. More precisely, we will design perturbations in four steps as follows.
\begin{itemize}
\item[$({\mathbf 1})$] Make the ends of $L\times I$ standard.
\item[$({\mathbf 2})$] Make $|x_0|$ non-increasing along any flow line of a positive function difference of $L\times I$ and decreasing outside a neighborhood of the Reeb chord manifolds and outside a neighborhood of $\{x_0=0\}$.
\item[$({\mathbf 3})$] Make the Reeb chords isolated in such a way that there are three Reeb chords of $L\times I$ corresponding to each Reeb chord manifold, one at $x_0=\pm 1$ and one near $x_0=0$.
\item[$({\mathbf 4})$] Make $L\times I$ generic with respect to rigid flow trees.
\end{itemize}
The goal of this perturbation process is to obtain a description of the rigid trees needed to compute the contact homology of $L\times I$ in terms of perturbed flow trees of $L\times I$.

Steps $({\mathbf 1})$ -- $({\mathbf 3})$ are straightforward to describe and are the subject of Subsection \ref{s:ends+Rch}. Step $({\mathbf 4})$ is more involved. It uses generalized flow trees, introduced in Subsection \ref{s:genfltr}, and will be completed in Subsection \ref{s:fltrgen}.

\subsection{Standard ends and isolated Reeb chords}\label{s:ends+Rch}
Consider Step $({\mathbf 1})$. Let $L\subset J^1(\R)$ be given by
$$
a(q)=\bigl(x_1(q),y_1(q),z(q)\bigr),\quad\quad q\in L.
$$
In coordinates $(x_0,y_0,x_1,y_1,z)$ on $J^1(\R^2)=T^\ast \R^2\times\R$, the unperturbed trace of the constant isotopy is then given by
$$
A(q,t)=\bigl(t,0,x_1(q),y_1(q),z(q)\bigr),\quad(q,t)\in L\times I.
$$
Fix a smooth even function $\psi\colon[-1,1]\to\R$ with the following properties: $\psi(0)=1$, $\psi$ and all its derivatives vanishes in $[0,\frac12]$, $\psi$ is decreasing on $[\frac12,1]$, and $\psi$ has a non-degenerate minimum at $1$, $\psi(1)=\frac12$. Let $\dot\psi$ denote the derivative of $\psi$. For $\eta>0$, define the Legendrian embedding $A^{({\mathbf 1})}_\eta\colon L\times I\to J^1(\R^2)$ by
\begin{equation}\label{e:A1etadef}
A^{({\mathbf 1})}_\eta(q,t)= \bigl(t,\eta\dot\psi(t)z(q),x_1(q),\eta\psi(t)y_1(q),\eta\psi(t)z(q)\bigr).
\end{equation}
The Reeb chords of $A^{({\mathbf 1})}_\eta(L\times I)$ are of two types: isolated Reeb chords in $\{x_0=\pm 1\}$ at the locations of the Reeb chords of $L$, and non-isolated Reeb chords along the manifolds $\hat M_c=M_c\cap\left\{|x_0|\le \frac12\right\}$.

Consider Step $({\mathbf 2})$. Fix a function $\phi\colon[-1,1]\to(-\infty,0]$, with the following properties: $\phi$ is even, $\phi(0)=0$, $\phi$ is decreasing on $(0,\frac12)$, $\phi$ and all its derivatives vanish at $[\frac12,1]$. Consider a neighborhood $N(M_c^\pm)$ of a Reeb chord endpoint manifold $M_c^\pm$ which, under $\Pi\colon J^{1}(\R^{2})\to\R^{2}$ maps to the product $\Pi(M_c)\times K$, where $K=\{x_1\colon |x_1-x_1(c)|<\delta\}$, where $x_1(c)$ is the $x_1$-coordinate of the Reeb chord $c$, for some small $\delta>0$. Consider a cut-off function $\hat \xi_K\colon K\to [0,1]$ which equals $0$ on $|x_1-x_1(c)|\le\frac12\delta$ and equals $1$ on $\frac34\delta<|x_1-x_1(c)|<\delta$. Let $\hat\xi\colon M_c\times K\to[0,1]$ denote the function with $\hat\xi|\{p\}\times K=\hat\xi_K$ for every $p$. Consider the pull-back of $\hat\xi$ to $N(M_c)$ and note that it extends constantly to all of $L\times I$. Denote the extension $\xi$. Write
$$
A^{({\mathbf 1})}_\eta(q,t)= \bigl(t,y_0^\eta(q,t),x_1^\eta(q,t),y_1^\eta(q,t),z^\eta(q,t)\bigr),
$$
where $A^{({\mathbf 1})}_\eta(q,t)$ is as in \eqref{e:A1etadef}.
The function $\phi$ depends only on the $x_0$-coordinate and the function $\xi$ only on the $x_1$ coordinate. Denote their respective derivatives by $\dot\phi$ and $\dot\xi$. Define the Legendrian embedding $A_\eta^{({\mathbf 2})}\colon L\times I\to J^1(\R^2)$ by
\begin{equation}\label{e:A2etadef}
A^{({\mathbf 2})}_\eta(q,t)=\bigl(t,(1+\eta\xi\phi)y_0^\eta+\eta\xi\dot\phi,
x_1^\eta,
(1+\eta\xi\phi)y_1^\eta+\eta\dot\xi\phi,
(1+\eta\xi\phi)z^\eta\bigr).
\end{equation}
If $\gamma$ is a flow line of a Legendrian submanifold $H\subset J^1(M)$ write $\gamma^+$ and $\gamma^-$ for the $1$-jet lift of $\gamma$ with the larger- and smaller $z$-coordinate, respectively.
\begin{lma}
Let $\eta>0$ and let $p$ be a point on a flow line $\gamma$ of a positive function difference of $A^{({\mathbf 2})}_\eta(L\times I)$. Assume that $p$ does not lie in $\{x_0=0\}$ and that not both $\gamma^+$ and $\gamma^-$ lie in the region where $\xi=0$. Then $|x_0|$ is decreasing along $\gamma$ at $p$.
\end{lma}

\begin{pf}
This follows from the choice of $\phi$ together with the fact that the $x_0$-component of $\nabla\xi$ is everywhere $0$.
\end{pf}

Consider Step $({\mathbf 3})$. In order to make Reeb chords of $L\times I$ isolated we fix functions $\beta_c\colon M_c\to(-\infty,0]$ for each Reeb chord $c$ of $L$. Consider $M_c=[-1,1]$ after projection to the $x_0$-line. Choose $\beta_c$ so that it has non-degenerate minima at $\pm 1$, a non-degenerate maximum with value $0$ close to $x_1=0$, and no other critical points. We assume moreover that the $x_0$-coordinates of the maxima of $\beta_c$ and $\beta_{c'}$ are different if $c\ne c'$. Furthermore let $\xi$ be a function just like $\xi$ discussed above but replacing $\delta$ by $\frac{\delta}{100}$ and let $\alpha=1-\xi$. We pull back the functions $\beta_c$ to neighborhoods $N(M_c^\pm)$ of the form describe above. Cutting these pull backs off with $\alpha$, we find that the function $\alpha\sum_c\beta_c$ extends in an obvious way to all of $L\times I$. We denote the extension $\alpha\beta$. Similarly, the functions $\dot\alpha\beta=\dot\alpha\sum_c\beta_c$ and $\alpha\dot\beta=\alpha\sum_c\dot\beta_c$ may be considered as functions on all of $L\times I$. Write
$$
A^{({\mathbf 2})}_\eta(q,t)= \bigl(t,y_0^\eta(q,t),x_1^\eta(q,t),y_1^\eta(q,t),z^\eta(q,t)\bigr),
$$
where $A^{({\mathbf 2})}_\eta(q,t)$ is as in \eqref{e:A2etadef}.
Define $A^{({\mathbf 3})}_\eta\colon L\times I\to J^1(\R^2)$ by
$$
A^{({\mathbf 3})}_\eta(q,t)
=\bigl(t,(1+\eta^2\alpha\beta)y_0^\eta+\eta^2\alpha\dot\beta,
x_1^\eta,
(1+\eta^2\alpha\beta)y_1^\eta+\eta^2\dot\alpha\beta,
(1+\eta^2\alpha\beta)z^\eta\bigr).
$$
Then $A^{({\mathbf 3})}_\eta(q,t)$ has isolated Reeb chords as claimed in $({\mathbf 3})$ above.

\subsection{Rough and fine scales -- limits on the rough scale}
As mentioned, we will consider the details of Step $({\mathbf 4})$, where we introduce our final perturbation of $L\times I$ resulting in a Legendrian embedding $A^{({\mathbf 4})}_\eta\colon L\times I\to J^{1}(\R^{2})$, in Subsection \ref{s:fltrgen}. In this subsection we will use rough properties of the perturbation of $A^{({\mathbf 4})}_\eta$ in order to derive rough results. More precisely, the property we use is the following.
\begin{itemize}
\item $A^{({\mathbf 4})}_\eta$ is a perturbation of $A_\eta^{({\mathbf 3})}$ of size $\Ordo(\eta^2)$ which is supported in the region outside the fixed neighborhoods of the Reeb chord manifolds.
\end{itemize}
In order to describe rigid flow trees of $A^{({\mathbf 4})}_\eta(L\times I)$, we consider convergence of such trees on two scales. On the rougher scale the important part of any flow trees concentrates in $x_0$-slices. This is fairly independent of the details of the perturbations in Step $({\mathbf 4})$. However, around such a slice of concentration we then re-scale the $x_0$-coordinate by $\eta^{-1}$ and study the corresponding microscopic limit as well. It is on this fine scale that the details of the perturbation in Step $({\mathbf 4})$ manifest themselves.

\begin{lma}\label{l:flinelimit}
Let $\gamma_\eta$, $\eta\to 0$ be a sequence of flow lines of $A_\eta^{({\mathbf 4})}(L\times I)$. Then, as $\eta\to 0$, some subsequence of $\gamma_\eta$ converges to a broken flow line on $L\times I$ consisting of flow lines of $L$ in horizontal slices and vertical flow lines of $\beta_c$ along $M_c$.
\end{lma}

\begin{pf}
Outside a neighborhood of $M_c$ the flow lines of $A^{({\mathbf 4})}_\eta(L\times I)$ obviously converge to horizontal curves which are parallel copies of flow lines of $L$. Inside a neighborhood of the manifolds $M_c$ the flow lines are as in a standard Morse-Bott situation, and any sequence of such flow lines converges to a combination of vertical flow lines of $\beta_c$ and horizontal flow lines.
\end{pf}

Lemma \ref{l:flinelimit} gives the local convergence on the rough scale. The following lemma describes the rough limits of flow trees of $A_\eta^{({\mathbf 4})}(L\times I)$ as $\eta\to 0$.

\begin{lma}\label{l:ftreelimit1}
Any sequence of flow trees, of $A_\eta^{({\mathbf 4})}(L\times I)$ of uniformly bounded formal dimension and which are not contained in $\{x_0=\pm 1\}$, has a subsequence that converges to a generalized flow tree on $L\times I$.
\end{lma}

\begin{pf}
Note that any gradient difference along a cusp edge of $A_\eta^{({\mathbf 4})}(L\times I)$ is transverse to the cusp edge. In particular, the preliminary transversality condition is met, see \cite[Subsection 3.1]{E}. It is a consequence of \cite[Lemma 3.12]{E}  and the preliminary transversality condition that the number of edges and vertices of trees in such a sequence is uniformly bounded. Hence by passing to a subsequence we may assume that the topological type of the trees in the sequence remains constant. The lemma then follows by applying Lemma \ref{l:flinelimit} to the edges of the trees in the sequence.
\end{pf}

As a first step toward the description of rigid flow trees of $A_\eta^{({\mathbf 4})}(L\times I)$ we rule out some generalized flow trees as limits of sequences of rigid trees.
\begin{lma}
Consider a sequence $\Gamma_\eta$ of rigid flow trees (with fixed punctures) of $A_\eta^{({\mathbf 4})}(L\times I)$ which converges to a generalized flow tree $G$. Let $\Gamma_t^r$, $r=1,\dots,m$ be the slice trees of $G$ and let $\hat n(G)$ denote the number of negative punctures of $G$ (or of $\Gamma_\eta$) which lies at $\hat c\in M_c$ for some $c$. Then
\begin{equation}\label{e:rest}
\sum_{r=1}^m(\dim(\Gamma_t^r)+1)-\hat n(G) = 0.
\end{equation}
\end{lma}

\begin{pf}
It follows from the dimension formula for flow trees that the formal dimension of a flow tree close to $G$ is given by the left hand side of \eqref{e:rest}.
\end{pf}
Note in particular that if the Legendrian knot $L$ is generic then $\dim(\Gamma^r_t)\ge 0$ for all trees in \eqref{e:rest}.

\subsection{Flow tree genericity}\label{s:fltrgen}
In this subsection we present the details of the construction of  $A^{({\mathbf 4})}_\eta\colon L\times I\to J^1(\R^2)$. The Legendrian submanifold $A^{({\mathbf 4})}_\eta(L\times I)$ will be generic with respect to rigid flow trees for all sufficiently small $\eta>0$ and its rigid flow trees will admit a description in terms of generalized flow trees with auxiliary data. In fact, we will first perturb $A^{({\mathbf 2})}_\eta$ outside a neighborhood of the Reeb chord manifolds in such a way that $A^{({\mathbf 4})}_\eta$ is obtained from the perturbed $A^{(\mathbf 2)}_\eta$ exactly as $A^{({\mathbf 3})}_\eta$ was obtained from the original $A^{({\mathbf 2})}_\eta$ in Subsection \ref{s:ends+Rch}.

To this end, we first describe the space of (partial) generalized trees near a slice tree. Let $a, b_1,\dots,b_k$ be Reeb chords of $L$ and consider $[\tau-\delta,\tau+\delta]\subset[-1,1]$, where $\tau=x_0(\hat b)$ for $b=b_j$ some $j$. Then by our choice of functions $\beta_c\colon M_c\to\R$, the gradient of $\beta_c$ induces an orientation on $M_c\cap\{|x_0-\tau|\le\delta\}$ for every $c\ne b$ provided $\delta>0$ is sufficiently small. Consider the compactified moduli space of trees
$$
{\overline \M}={\overline \M}(a; b_1,\dots,b_k)
$$
of $L$ with positive puncture at $a$ and negative punctures at $b_1,\dots,b_k$. The interior of this space is a manifold of dimension $d$. (Local coordinates can be obtained by the location of the branch points on the boundary of the corresponding holomorphic disks.) The boundary $\pa{\overline \M}$ of ${\overline \M}$ consists of products of lower-dimensional moduli spaces corresponding to broken disks. If $\M_1\times \M_2\times\dots\times \M_r$ is such a product then the positive puncture of $\M_1$ is $a$, the positive puncture of $\M_j$ agrees with a unique negative puncture of some $\M_k$, for $k<j$, the negative punctures which are not matched by some positive puncture are $b_1,\dots,b_k$, and the dimensions $d_j$ of $\M_j$ satisfy
\begin{equation}\label{e:dimglue}
\sum_{j=1}^r d_j + (r-1) = d.
\end{equation}
As with generalized flow trees it is convenient to organize such a product of moduli spaces into levels, where $\M_1$ is the first level and all moduli spaces with positive puncture paired with a negative puncture of $\M_1$ constitute the second level, and, in general, all moduli spaces with positive puncture paired with a negative puncture of a moduli space of level $j$  constitute the $(j+1)^{\rm th}$ level. In this way ${\overline \M}$ is a naturally stratified space (a manifold with boundary with corners endowed with a Kuranishi structure, see \cite{FOOO} for this notion). We will associate a $(d+1)$-dimensional stratified space of (partial) generalized flow trees to $\M$.  We call this space $\Pe=\Pe(a;b_1,\dots,b_k)$. It is constructed from strata of ${\overline \M}$ exactly as the space of generalized flow trees except that gradient lines need not be complete, as follows. Consider first the top-dimensional stratum. We associate to this stratum the space $\Pe_0$ which consists of trees $\Gamma_t$ in $\M(a;b_1,\dots,b_k)$ viewed as a tree in a slice, $t\in[\tau-\delta,\tau+\delta]$. This is a $(d+1)$-dimensional space: $d$ dimensions for the tree and one dimension for the slice. Consider next a stratum $S$ with $k$ levels $L_1,\dots,L_k$. We associate a $(d+1)$-dimensional space $\Pe_S$ of partial generalized trees to this stratum as follows. First consider the tree in $L_1$ as a slice tree $\Gamma_t$. For each negative puncture of $\Gamma_t$ which matches a positive puncture of some tree in $\M_2$, pick an oriented flow segment of the corresponding Reeb chord manifold and let the tree in the second level have its positive puncture where this flow line ends and continue in this way until $\M_k$ is reached. It follows from \eqref{e:dimglue} that the dimension of $\Pe_S$ equals $(d+1)$. Moreover, the spaces constructed fits together in an obvious way to a $(d+1)$-dimensional space which is our space $\Pe$. In a sense it is a resolution of ${\overline \M}$.

Our next objective is to define an evaluation map
$$
\ev^F\colon \Pe(a;b_1,\dots,b_k)\to M_{b_1}\times\dots\times M_{b_k},
$$
where $F$ labels a specific perturbation of the Legendrian $L\times I$. This map is constructed inductively.

The starting point for the construction is the observation that it is easy to reconstruct actual flow trees of $A^{({\mathbf 2})}_\eta(L\times I)$, see Subsection \ref{s:ends+Rch}, near flow trees in slices. More precisely, let $F\colon L\times I\to J^1(\R^2)$ be a deformation of $L\times I$. We are only interested in small perturbations so let its size be $\Ordo(\eta)$ as $\eta\to 0$. Assume further that the perturbation is supported in the complement of the Reeb chord manifolds $M_c$ (just like for $A^{({\mathbf 2})}_\eta$ in Subsection \ref{s:ends+Rch}). Outside the neighborhood of the Reeb chord manifolds, the $x_1$-component of the gradient of any local function difference of $L\times I$ is bounded from below. It follows from this that the projection of any flow tree on $F(L\times I)$ for $\eta>0$ sufficiently small, to a slice is a tree.

Conversely, given a tree $\Gamma'$ in a slice, we can integrate it to a tree on $F(L\times I)$ for small $\eta$. The integration procedure is inductive and amounts to solving the gradient equation corresponding to the perturbation $F$ along the edges of the tree. For the $x_1$-component of the edge this results simply in a re-parametrization and the $x_0$-component becomes non-constant along the edge. To extend the integration over the entire tree we start integrating along the edge (edges) emanating from the positive puncture of the tree with the initial condition given by the $x_0$-coordinate of the slice of $\Gamma'$ and follow the flow orientation. Inductively, the result of integrations along edges closer to the positive puncture in the tree give initial conditions for the integrations along outgoing edges at any vertex in the tree. If $\Gamma'\in\M(a;b_1,\dots,b_k)$, if the punctures of $\Gamma$ are $p_0,p_1,\dots,p_k$ with $p_0$ the positive puncture, and if $\Gamma'$ is considered a tree in a slice with $x_0$-coordinate $\alpha$, then we write $\iota(p_0, p_j)\in M_{b_j}$ for the result of integration along $\Gamma'$ at the negative puncture $p_j$ with initial condition given by the $x_0$-coordinate at the positive puncture $p_0$.

Let $\Gamma\in\Pe(a;b_1,\dots,b_k)$ be a partial generalized flow tree with levels $(\Gamma_j^{1},\dots,\Gamma_j^{s_j})$, $j=1,\dots,r$. We define $\ev^{F}(\Gamma)$ as follows. Apply the integration procedure to the top level $\Gamma_1^{1}$, starting at the $x_0$-coordinate of its positive puncture in $M_{a}$. This gives points $\iota(p_0^{1},p_j^{1})\in M_{c_j}$ for all the negative punctures $p_j^{1}$ of $\Gamma_1^{1}$. If a negative puncture $p_j^{1}$ of $\Gamma^{1}_1$ is also a negative puncture of $\Gamma$ then define $\iota'(p_j^{1})=\iota(p_j^{1})\in M_{c_j}$ to be the corresponding $M_{c_j}$-component of the evaluation map $\ev^{F}$. If a negative puncture $p_k^{1}$ of $\Gamma_1^{1}$ is not a negative puncture of $\Gamma$ then $\iota(p_0^{1},p_k^{1})\in M_{c_k}$ is the initial condition for the integration procedure of the second level tree attached at $p_k^{1}$. Continuing inductively in this way over all the levels we get for each negative puncture $p_j$ of $\Gamma$, which is a negative puncture $q_r$ of some level tree in $\Gamma$ with positive puncture $q_0$, a point $\iota'(p_j)=\iota(q_0,q_s)\in M_{b_j}$. Define
\[
\ev^{F}(\Gamma)=\left(\iota'(p_1),\dots,\iota'(p_k)\right)\in M_{b_1}\times\dots\times M_{b_k}.
\]

Let $b$ be a Reeb chord, let $\tau=x_0(\hat b)$, and let $p_0,p_1,\dots,p_k$ denote the punctures of $\Gamma\in\Pe(a;b_1,\dots,b_k)$. Let $J=\{p_{j_1},\dots,p_{j_s}\}$ be a subset of the punctures of $\Gamma$. Consider the subvariety $\Sigma_{J,\tau}\subset M_{b_1}\times\dots\times M_{b_k}$,
$$
\Sigma_{J,\tau}=\left\{(t_1,\dots,t_k)\colon t_j=\tau,
\text{ for all }j\text{ such that }p_j\in J\right\}.
$$
We show next that there are small perturbations $F$ of ${A^{({\mathbf 2})}_\eta}$ such that $\ev^{F}$ is transverse to $\Sigma_{J,\tau}$ for all small $\eta>0$. We will write $\ev^{F}_J$ for the evaluation map $\ev^{F}$ composed with the projection to $\Pi_{p_j\in J}M_{b_j}$ (i.e. to the product of the factors corresponding to punctures in $J$).

We start with the following preliminary lemma concerning a general fact about flow trees with only one positive puncture.

\begin{lma}\label{l:onlyonce}
Let $\Gamma$ be a flow tree of a Legendrian submanifold $L\subset J^1(M)$ with only one positive puncture and let $v$ be a vertex of $\Gamma$. Then any sheet of $L$ which contains some point $\tilde v$ in the $1$-jet lift of $v$ contains the $1$-jet lift of exactly two flow lines with endpoint at $\tilde v$.
\end{lma}

\begin{pf}
We label a sheet by its corresponding local function. As mentioned above, each edge in a tree with one positive puncture is naturally oriented by using the negative gradient flow of the local function difference which is positive and with this convention exactly one edge at each vertex is oriented toward it (incoming flow line), all other away from it (outgoing flow lines).

Assume that some sheet, $f_1$ say, contains lifts of more than two flow lines. As the $1$--jet lift is an oriented curve there must exist at least two lifts which are oriented away from $\tilde v$. By the definition of a flow tree there must exist matching edges oriented toward $\tilde v$ and at most one of these can be a lift of the incoming edge. The remaining one must come from an outgoing edge which thus is a flow between $f_2$ and $f_1$ where $f_2>f_1$. Continuing this argument with the matching lift of the $1$-jet lift segment in $f_2$ etc, we find that the tree has at least two incoming edges at $v$ since there are only finitely many sheets. This contradiction establishes the lemma.
\end{pf}

Next we derive certain injectivity properties of flow trees of $A^{({\mathbf 2})}_\eta(L\times I)$.
\begin{lma}\label{l:discset}
Assume that the set of points where two local gradient differences of $L$ agree is discrete (this is obviously an open condition on $L$). If $v$ is any vertex of a flow tree $\Gamma$ of $A^{({\mathbf 2})}_\eta(L\times I)$ with incoming edge $e^i$ and outgoing edges $e^o_1,\dots,e^o_m$ then the self intersections of the $1$-jet lift of $\Gamma$ near the $1$-jet lift of $v$ form a discrete set.
\end{lma}

\begin{pf}
Consider an arc in one of the sheets. It is a consequence of Lemma \ref{l:onlyonce} that there is one incoming and one outgoing flow in this sheet. The genericity condition on $L$ implies that the flows of the corresponding gradients intersect only at a discrete set of points.
\end{pf}

In order to establish the desired transversality, we will work with a re-scaled version of the map $\ev_J^{F}$,
$$
s_\eta[\ev_J^{F}]=\frac{1}{\eta}\left(\ev_J^F-(\tau,\dots,\tau)\right)\colon \Pe\to\R^{|J|},
$$
where $|J|$ denotes the cardinality of $J$. Here we will take $\eta\to 0$, and apply the Sard-Smale theorem. Let $\Fu$ be the space of perturbations of $L\times I$ as described. (To get a hold of this space we can think of it as of the space of all functions on $L\times I$ supported outside a neighborhood of the $M_c$.) The evaluation map discussed above then gives a map
$$
\ev_J\colon \Fu\times \Pe\to M_b\times\stackrel{r}{\cdots}\times M_b,
$$
with a re-scaling $s_\eta[\ev_J]$ analogous to that defined above.

\begin{lma}\label{l:tv}
The map $s_\eta[\ev_J]$ is transverse to $0\in\R^{|J|}$ in a neighborhood of $A^{({\mathbf 2})}_\eta$ for all sufficiently small $\eta>0$.
\end{lma}

\begin{pf}
Since transversality is an open condition we need only show that the differential is onto at any $(F,\Gamma)\in \Fu\times\Pe$ with $s_\eta[\ev_J^{F}](\Gamma)=0$. Fix such $(F,\Gamma)$. To prove surjectivity we show that $\pa_j+\Ordo(\eta)$ is in the image if the differential for each $j$, where $\pa_j$ is the tangent vector along the $j^{\rm th}$ factor in the product $(M_b)^{r}$.

Consider the $j^{\rm th}$ puncture $q_j$ mapping to $\tau\in M_b$. Assume for definiteness that $\tau<0$. Start at $q_j$ and follow $\Gamma$ in the direction opposite to the flow orientation. Then we meet a first edge $e$ which leaves the region $N(M_b)$ around $M_b$ where no perturbation is supported. Perturb this edge near the entrance point of $N(M_b)$ by deforming the Legendrian submanifold. We claim that if the region $N(M_b)$ is chosen small enough then the only negative puncture in $\Gamma$ which appears after the edge $e$ in the flow orientation and which map to $\tau\in M_b$ is $q_j$. To see this, we use the fact that the distance between any two negative punctures in a tree is uniformly bounded from below, see Lemma \ref{l:distpunct} below, to conclude that to get to a negative puncture after $q_j$, we must follow some edge which leaves $N(M_b)$ and in such a region, integration along any edge induces a strict decrease the $x_0$-coordinate. Thus, the $x_0$-coordinate is smaller than $\tau$ at any later puncture.  

Consider next punctures $q_r$ which appear after $q_j$, maps to $\tau\in M_b$ but lies on a lower level tree. For such punctures we can remove the corresponding components of $\pa_r$ which the already introduced perturbation gives rise to by varying the length of the attaching flow line.

Finally we must deal with punctures connected to $e$ only via vertices above $e$. (The perturbation could effect such a vertex.) Consider the vertices connected to $e$ by an edge attached as close to $e$ as possible. These vertices lie above the perturbation region in $e$. After Lemma
\ref{l:discset}, we know that $\Gamma$ has injective points in this region and it is therefore easy to compensate for the shifts induced by earlier perturbations. Using an obvious induction we take care of all punctures in the tree above $q_j$. (Note that all perturbations needed for these additional punctures appear off the main stem connecting $q_j$ to the positive puncture and at $x_0$-levels above that of $e$.) 
\end{pf}

\begin{cor}\label{c:mfd}
For an open dense set of perturbations $F$ of $A^{({\mathbf 2})}_\eta$,
$$
(\ev_J^F)^{-1}\left((\tau,\stackrel{r}{\dots},\tau)\right),
$$
where $\ev^{F}_J\colon\Pe(a;b_1,\dots,b_k)\to M_b\times\stackrel{r}{\dots}\times M_b$, $\dim(\Pe (a;b_1,\dots,b_k))=d+1$ is a manifold of dimension $d+1-r$.
\end{cor}
\begin{pf}
This is a standard application of the Sard-Smale theorem.
\end{pf}

We next specialize to the situation of most importance to us. Consider a compactified moduli space ${\overline \M}(a;b_1,\dots,b_k)$ of flow trees on $L$ of dimension $d$ and with at least $m\ge d$ negative punctures mapping to $b$ (i.e. $b_j=b$ for at least $d$ indices $j$). Consider the corresponding space $\Pe=\Pe(a; b_1,\dots,b_k)$ and let $J$ be a subset of the punctures of trees in $\Pe$ mapping to $M_b$ of cardinality $|J|=d$. Note that the map $s_\eta[\ev^{F}_J]\colon\Pe\to\R^{|J|}$ is continuous on the space $s_\eta[\Pe]$ obtained from $\Pe$ by re-scaling the lengths of all connector flow lines by $\eta^{-1}$. Moreover, it is clear that the closed subset $\left(s_\eta[\ev_J^F]\right)^{-1}(0)\subset s_\eta[\Pe]$, where $F$ is of size $\Ordo(\eta)$, is bounded. Hence it is compact. Thus, Corollary \ref{c:mfd} implies that, for generic $F$, $\left(\ev^{F}_J\right)^{-1}(\tau,\dots,\tau)$ is empty if $|J|>d+1$ and is a compact $0$-dimensional manifold if $|J|=d+1$. Furthermore, any point in such a $0$-dimensional manifold $\left(\ev^{F}_J\right)^{-1}(\tau,\dots,\tau)$ corresponds to a point in the interior $\M=\M(a;b_1,\dots,b_k)$ of ${\overline \M}={\overline \M}(a; b_1,\dots,b_k)$ for the following reason. Recall that the boundary of the $d$ dimensional ${\overline \M}$ consist of products of the form
$$
\M_{d_1}\times\cdots\times\M_{d_r},
$$
where $\M_{d_j}$ is a moduli space of dimension $d_j\ge 0$ and where
$$
d_1+\cdots+d_r+(r-1)=d.
$$
Thus if $\left(\ev^{F}_J\right)^{-1}(\tau,\dots,\tau)$ intersects $\pa{\overline \M}$ then there are subsets $J_l$, $J=\cup_l J_l$ of negative punctures of trees in $\M_{d_l}$ for which $\ev^{F}_{J_l}(\tau,\dots,\tau)$ is non-empty. The genericity condition then implies that
$|J_l|\le d_l+1$. Then $d-(r-1)=\sum_l d_l\ge \sum_l|J_l|-1=|J|-1$ and $d+1>|J|$ since $r\ge 2$. This however contradicts $|J|=d+1$ and thus shows that all points lie in $\M$.

Let $p_0,p_1,\dots,p_k$ denote the punctures of the generalized flow trees in $\Pe(a; b_1,\dots,b_k)$, where $\dim(\M(a;b_1,\dots,b_k))=d$, let $F$ be a generic perturbation, let $J$ denote a subset of the punctures mapping to the Reeb chord $b$ with $|J|=d+1$, and let $x_0(\hat b)=\tau$. We write $E_{J}^{F}\left({\overline \M}(a;b_1,\dots,b_k)\right)=\left(\ev^{F}_J\right)^{-1}(\tau,\dots,\tau)$.

\begin{rmk}
It is easy to see that if $\dim(\M(a;b_1,\dots,b_k))=0$ then $E^{F}_J({\overline \M}(a;b_1,\dots,b_k))$ consists of exactly one point.
\end{rmk}

Since the moduli space of flow trees of $L$ is compact and since there are only finitely many Reeb chord manifolds, Lemma \ref{l:tv} implies that we can find a perturbation $F$  of $A_\eta^{({\mathbf 2})}$ which makes all evaluation maps $\ev^{F}\colon{\overline \M}(a;b_1,\dots,b_k)$ transverse to all varieties $\Sigma_{J,\tau}$, $\tau=x_0(b)$ some $b$ and $J$ any subset of the negative punctures mapping to $b$, for all moduli spaces and all choices of collections of punctures.

We fix a choice of such an $F$ and perform also the perturbations of size $\Ordo(\eta^2)$ along the Reeb chord manifolds which were used to obtain $A^{({\mathbf 3})}$ from $A^{({\mathbf 2})}$ in Subsection \ref{s:ends+Rch}. Denote the resulting Legendrian embedding $A^{({\mathbf 4})}_\eta\colon L\times I\to J^1(\R^2)$. We next give a description of the rigid flow trees of $A^{({\mathbf 4})}_\eta(L\times I)$ for small $\eta>0$.

Consider a generalized flow tree $G\in\Pe(a;b_1,\dots,b_s)$ on $L\times I$. A {\em complete slice} $S$ of $G$ is a collection of slice trees of $G$ of the following form. 
\begin{itemize}
\item
There is one slice tree $\Gamma^{0}$ in $S$ which lies on a level above (i.e. a level of lower numbering than) all other slice trees in $S$. 
\item
The positive puncture of $\Gamma^{0}$ is either the positive puncture of $G$ or connected to a negative puncture of some other slice tree in $G$ by a flow line of length $>0$. 
\item
Every slice tree in $S$ of level $k$ which is lower than the level of the top slice-tree $\Gamma^{0}$ is connected at its positive puncture {\em via a flow line of length $0$} to the negative puncture of some slice tree in $S$ of level $k-1$. 
\end{itemize}
Thus, a complete slice $S$ is a (possibly broken) tree in a level lying in some compactified moduli space ${\overline \M}_S={\overline \M}(e;c_1,\dots,c_k)$. If the complete slice $S$ lies in the $x_0$-slice of $\hat c$ where $c=c_j$ for one of the negative punctures $c_j$, then $E_{J}^{F}({\overline \M}_S)$ may be non-empty for suitable $J$.

Consider a generalized flow tree $G\in\Pe(a;b_1,\dots,b_s)$ of $L\times I$ and let $S_1^G,\dots,S_k^G$ be its complete slices. We say that $G$ is {\em potentially rigid} if for every $S_k^{G}\in{\overline \M}(e;c_1,\dots,c_k)$ in an $x_0(\hat c)$-slice, $c=c_j$ some $j$, the number of punctures mapping to $c$ (i.e. the number of indices $j$ such that $c=c_j$) is larger than $\dim(\M_{S_j^G})$. In this case we take $E^{F}({\overline \M}_{S_j^G})$ to be the union of all the manifolds $E_{J}^{F}({\overline \M}_{S_j^G})$ over all distinct choices of subsets $J$ of negative punctures with $|I|=d+1$.

\begin{lma}\label{l:incspace}
For all sufficiently small $\eta>0$, there is a $1-1$ correspondence between, on the one hand, rigid flow trees of $A^{(\mathbf 4)}_\eta(L\times I)$ and, on the other, the union over all potentially rigid generalized flow trees $G$ of the product sets
\begin{equation}\label{e:incspace}
E^{F}({\overline \M}_{S_1^G})\times E^{F}({\overline \M}_{S_2^G})\times\dots\times E^{F}({\overline \M}_{S_k^G}),
\end{equation}
where $S_j^G$, $j=1,\dots,k$, are the complete slices of $G$.
\end{lma}

\begin{pf}
Let $G$ be a potentially rigid generalized flow tree and let $F$ denote the perturbation of $A^{({\mathbf 2})}_\eta$ discussed above. A point in the product
$$
E^{F}({\overline \M}_{S_1^G})\times E^{F}({\overline \M}_{S_2^G})\times\dots\times E^{F}({\overline \M}_{S_k^G}),
$$
corresponds to $k$ complete slices of $G$ each of which give a transverse solution to an equation $s_\eta[\ev_{J}^{F}]=0$ for some Reeb chord $b$ and some collection of punctures $J$. The Legendrian embedding $F$ is at $\Ordo(\eta)$ distance from the inclusion of $L\times I$ and has Morse-Bott Reeb chords along $\hat M_c$ for all Reeb chords $c$ of $L$. To obtain $A^{({\mathbf 4})}_\eta$ from $F$ we make a perturbation of size $\Ordo(\eta^2)$ in a small neighborhood of the manifolds $\hat M_c$.

The evaluation map is re-scaled by $\eta^{-1}$. After this re-scaling, the perturbation from $F$ to $A^{({\mathbf 4})}_\eta$ is a standard perturbation out of a Morse-Bott situation into a Morse situation. In particular, flow lines before the perturbation with evaluation map transversely equal to the location of the maximum corresponds in a $1-1$ fashion to flow lines after the perturbation ending at the created maximum. We conclude that for $\eta>0$ sufficiently close to $0$, any point in the product gives rise to a collection of partial flow trees with one negative puncture at a $\hat c$ for each negative puncture which maps to $x_0(\hat c)$ by the evaluation map.

We show that these pieces can be glued in a unique way to a rigid flow tree of $A^{({\mathbf 4})}_\eta(L\times I)$. This is straightforward: by the definition of complete slice, all complete slices are connected to other complete slices by flow lines of length bounded from below. In order to glue the pieces, we thus need only make sure that the gluing problem for the flow lines of $\beta_c$ connecting the negative puncture at a complete slice of an above level to a complete slice in level below has a unique solution. Consider first the case when the critical point is a maximum. Since the Morse-Bott perturbation is of size $\Ordo(\eta^2)$, a change in the tree of the higher level near the negative puncture of size $\Ordo(\eta^2)$ produces a finite change in the level of the outgoing flow line. Since the solution to $s_\eta[\ev_c^{F}]=0$ on the higher level is uniformly transverse, it persists under changes of size $\Ordo(\eta^2)$. The case when the negative puncture has minimum character is similar: in order to connect the subtrees only perturbations of the pieces of size $\Ordo(\eta^2)$ are needed and the solutions to $s_\eta[\ev_b^{F}]=0$ persists under such perturbations. We conclude that near each point in the product \eqref{e:incspace} there is a unique rigid flow tree of $A^{({\mathbf 4})}_\eta(L\times I)$ for all $\eta>0$ sufficiently small.

In order to finish the proof it remains only to check that any sequence of rigid trees of $A^{({\mathbf 4})}_\eta(L\times I)$ converges to a potentially rigid generalized tree and gives a solution to the evaluation condition of \eqref{e:incspace}. The fact that it converges to a potentially rigid tree is an easy consequence of Lemma \ref{l:ftreelimit1}. Re-scaling by $\eta^{-1}$ takes us as mentioned above to a standard Morse-Bott situation near the Reeb chord manifolds and it follows that any sequence of trees give a solution.
\end{pf}

We call the product of $0$-dimensional manifolds as in \eqref{e:incspace} {\em incidence spaces}. Thus, the conclusion of Lemma \ref{l:incspace}, is that the contact homology differential of $A^{({\mathbf 4})}_\eta(L\times I)$, for $\eta>0$ sufficiently small, admits a description in terms of generalized flow trees and their corresponding incidence spaces. Although this is a rather nice geometric description, algebraically, it is in general messy. Moreover, even geometrically it is not very explicit: it is, in general, difficult to determine the incidence numbers since they require exact knowledge of moduli spaces of flow trees on $L$ of arbitrary dimension. (For a fixed knot $L$ this is less of a problem since any $L$ admits a presentation in which there are no disks with multiple negative punctures and for such knots the description of the algebra is straightforward.) However, we need to deal with one parameter families of knots, and want the algebra as simple as possible and as we shall see, there is an algebraically much preferable description of an algebra which is stable tame isomorphic to the one discussed above. This algebra arises from abstract perturbations and is related to the one above directly through its description in terms of generalized flow trees and incidence spaces. For this reason we note the following.

\begin{rmk}\label{r:geomasabstr}
Fix a small $\eta>0$ for which Lemma \ref{l:incspace} holds. Consider the Legendrian submanifold $F(L\times I)$ where $F$ is a small perturbation of $A^{({\mathbf 2})}_\eta$. The integration along trees discussed above map can be viewed as a map ${\overline \M}\times I\to C^k(S^1,T^{\ast} I)$,  where $S^1$ is naturally identified with the cotangent lift of the tree. Lemma \ref{l:incspace} then gives a description of the moduli space of rigid trees determined by $A^{({\mathbf 4})}_\eta(L\times I)$ in terms of perturbed trees and Morse functions $\beta_c\colon M_c\to\R$ on the Reeb chord manifolds, where the perturbation function is the $x_0$-component of the integration map discussed above.
\end{rmk}

\section{Abstract perturbations and contact homology computations}\label{S:abstrprt}
In this section we design a perturbation for the trace of the constant isotopy which leads to a simple differential.

\subsection{A perturbation function}
Let $\Gamma$ be a flow tree on $L\subset J^1(\R)$. Note that the cotangent lift of $\Gamma$ is naturally  by a map of a circle $S_\Gamma$ with metric induced by the parametrization and with ordered marked points $p_0, p_1,\dots,p_k$ at the punctures of $\Gamma$. (Here $p_0$ corresponds to the positive puncture.) We will associate to $\Gamma$ a function $\phi_\Gamma\colon S_\Gamma\to\R$ such that $\phi_\Gamma\ge 0$ and such that
$$
\phi_\Gamma(p_0)=0=\phi_\Gamma(p_k)<\cdots<\phi_\Gamma(p_1).
$$
Furthermore, the functions $\phi_\Gamma$ will vary continuously with $\Gamma$ in the moduli space of flow trees on $L$. In order to get a common source circle for all the maps $\phi_\Gamma$, we scale all cotangent lift circles $S_\Gamma$ so that they have length $2\pi$ and view $\Gamma\mapsto\phi_\Gamma$ as a map ${\overline \M}\to C^k(S^1,\R)$, where ${\overline \M}$ is the compactified moduli space of trees of $L$ and where $k\ge 1$. Before discussing these functions we prove a preliminary lemma.

\begin{lma}\label{l:distpunct}
There exists a constant $C>0$ such that for any flow tree $\Gamma$ of $L$ with one positive puncture, the distance in $S_\Gamma$ between any two punctures is bounded below by $C$.
\end{lma}

\begin{pf}
Let $c$ be a Reeb chord of $L$. The only Reeb chord in a neighborhood of $c$ is $c$ itself. For area reasons, if $c$ is the positive puncture of a tree then it cannot also appear as a negative puncture in
that tree. Thus it suffices to consider two neighboring negative punctures in a tree $\Gamma$ which both map to $c$. Consider the $1$-jet lift $\tilde\Gamma$ of $\Gamma$. Since the path of $\tilde\Gamma$ which connects the punctures is incoming at one of them and outgoing at the other, it connects the top endpoint of $c$ to the bottom endpoint of $c$ and thus its length is bounded from below.
\end{pf}

Assume that $L\subset J^1(\R)$ is sufficiently generic so that the lengths of all its Reeb chords are pairwise distinct. Let the lengths be
$$
0<l_1<l_2< \cdots< l_m.
$$
It follows by Stokes theorem that there exists $M$ (take $M>\frac{l_m}{l_1}$) such that no flow tree with one positive puncture has more than $M$ negative punctures, and that there exists $K$ (take $K>\frac{l_m}{\min_{j}\{l_{j+1}-l_j, l_1\}}$) such that no broken tree has more than $K$ levels.

Choose a function $h\colon\R\to\R$, $h>0$, such that for every $j=1,\dots,r$
$$
\tfrac{1}{2KM}h(l_{j+1})> h(l_j).
$$
We will define the map ${\overline \M}\to C^k(S^1,\R)$ in an inductive manner using the area filtration of ${\overline \M}$. In the first step we define the map on the moduli space of flow trees of smallest area. In later steps we define the map on the compactification of the moduli space of trees of higher areas assuming that the map is already defined for all moduli spaces of trees of smaller area. Since the compactification of a moduli space of trees of a given area consists of broken trees with pieces of smaller areas, the map has a natural definition on the boundary of the moduli space and we show how to extend it to the interior.

Let $\Gamma$ be a flow tree. We say that a smooth function $\phi_\Gamma\colon S_\Gamma\to\R$ is {\em stretching} if the following holds
\begin{itemize}
\item If $\Gamma$ has one or zero negative punctures, then $\phi_\Gamma=0$.
\item If $\Gamma$ has $r>1$ negative punctures and if $l$ is the length of its positive Reeb chord, then $\phi_\Gamma(p_0)=\phi_\Gamma(p_r)=0$, $\phi_\Gamma$ increases monotonically as we move along the circle in the negative direction from the point midway between $p_r$ and $p_{r-1}$ to the point $q$ midway between $p_1$ and $p_0$, $\phi_\Gamma$ decreases monotonically between $q$ and the point midway between $q$ and $p_0$ and then it is constant. Furthermore,
$\phi_\Gamma(p_{j})-\phi_\Gamma(p_{j+1})=\frac{1}{M}h(l)$ for $j=1,\dots,r-1$.
\end{itemize}

Let
\begin{equation}\label{e:slicedecomp}
\Gamma=\Gamma_1^1
\cup\left(\Gamma_2^1\cup\dots\cup\Gamma_2^{s_2}\right)
\cup\dots
\cup\left(\Gamma_l^1\cup\dots\cup\Gamma_l^{s_l}\right)
\end{equation}
be a broken tree where $\Gamma^{1}_j,\dots,\Gamma^t_j$ are the trees of level $j$. Let $\phi_{\Gamma_j^t}\colon S_{\Gamma_j^t}\to\R$ be smooth functions which are constant in neighborhoods of punctures. Then these functions glue in an obvious way to a smooth function $\phi_\Gamma\colon S_\Gamma\to\R$ as follows. For points $x\in\Gamma^{1}_{1}\subset\Gamma$, let $\phi_\Gamma(x)=\phi_{\Gamma^{1}_{1}}(x)$. Assume inductively that $\phi_\Gamma$ has been defined for all points in trees of level $<j$. Let $x\in\Gamma_j^t$ and let $q_{j}^{t}\in\Gamma_{j-1}^{r}$ denote the negative puncture in the level $j-1$ tree $\Gamma_{j-1}^{r}$ where the positive puncture of $\Gamma_j^{t}$ is attached. Define $\phi_\Gamma(x)=\phi_\Gamma(p_j^{t})+\phi_{\Gamma_{j}^{t}}(x)$.

We define the {\em common ancestor $A_\Gamma(q,q')$} of two negative punctures $q,q'$ in the broken tree $\Gamma$ as an unbroken subtree of $\Gamma$, inductively, in the following way. If both $q$ and $q'$ are negative punctures of the same unbroken tree $\Gamma^t_j$ then we take $A(q,q')=\Gamma^t_j$. If $q$ and $q'$ lie on different levels $l(q)$ and $l(q')$ with $l(q)<l(q')$, say, then we take $A_\Gamma(q,q')=A_{\Gamma'}(q,q'')$ where $q''$ is the negative puncture of the tree of level $l(q')-1$ where the positive puncture of the tree of $q'$ is attached and where $\Gamma'$ is the broken tree obtained from $\Gamma$ by cutting at $q''$. If $q$ and $q'$ lie on the same level but not in the same tree then we take $A_\Gamma(q,q')=A_{\Gamma'}(p,p')$, where $p$ and $p'$ are the negative punctures in the trees of level $l(q)-1=l(q')-1$ where the positive punctures of the trees of $q$ and $q'$, respectively, are attached, and where $\Gamma'$ is the broken tree obtained from $\Gamma$ by cutting at $p$ and at $p'$.

Let $p_0,p_1,\dots, p_r$ denote the punctures of a broken tree $\Gamma$ as in \eqref{e:slicedecomp} .
\begin{lma}\label{l:jumpbroken}
If $\phi_{\Gamma_j^{t}}$ are stretching for all $j$ then for any $t$, $1\le t\le r-1$
$$
\phi_\Gamma(p_{t})-\phi_\Gamma(p_{t+1})>\left(1-\frac{L}{(K+1)}\right)\frac{1}{M}h(l(p_t,p_{t+1})),
$$
where $L$ is the number of levels in $\Gamma$ and where $l(p_t,p_{t+1})$ is the length of the Reeb chord of the positive puncture in $A_\Gamma(p_t,p_{t+1})$.
\end{lma}

\begin{pf}
We use induction. For trees of only one level this is immediate. Assume next that it holds for all broken trees of $L-1$ levels and consider attaching an $L^{\rm th}$ level to an $(L-1)$-level broken tree $\Gamma'$ to form a $L$-level tree $\Gamma$. The only point that needs to be checked is that if $q$ is the first negative puncture in a tree $\Delta$ of level $L$ and if $p$ is the negative puncture in $\Gamma$ preceding it in $\Gamma$, then
$$
\phi_\Gamma(p)-\phi_\Gamma(q)>\left(1-\frac{L}{K+1}\right)\frac{1}{M}h(l(p,q)).
$$
Note that either $p=p'$ where $p'$ lies in a tree of level at most $L-1$ or $p$ lies in a tree of level $L$ attached at a puncture $p'$ of some tree of level $L-1$. In the latter case we have $\phi_\Gamma(p)=\phi_{\Gamma'}(p')$. If $q'$ denotes the negative puncture in the tree $\Gamma'$ at which the tree $\Delta$ is attached then the inductive assumption implies
$$
\phi_{\Gamma'}(p')-\phi_{\Gamma'}(q')
<\left(1-\frac{L-1}{K+1}\right)\frac{1}{M}h(l(p',q')).
$$
By the definition of common ancestor $A(p',q')=A(p,q)$ and hence $h(l(p',q'))=h(l(p,q))$. Thus, since $\phi_\Delta$ is stretching, with $l=l(p,q)$,
\begin{align*}
\phi_\Gamma(p)-\phi_{\Gamma}(q) &= \phi_{\Gamma'}(p')-\phi_{\Gamma'}(q')-\phi_\Delta(q)\\
&>\left(1-\frac{L-1}{K+1}\right)\frac{1}{M}h(l)-M\frac{1}{2KM}h(l)\\
&>\left(1-\frac{L-1}{K+1}-\frac{1}{K+1}\right)\frac{1}{M}h(l)=
\left(1-\frac{L}{K+1}\right)\frac{1}{M}h(l).
\end{align*}
\end{pf}

\begin{lma}\label{l:shiftfct}
There exists a continuous function $\M\to C^k(S^1,\R)$, $\Gamma\mapsto\phi_\Gamma$, such that the following holds. If $\Gamma$ is any (broken) tree in $\M$ with negative punctures $p_1,\dots,p_r$ then
$$
\phi_\Gamma(p_{t})-\phi_\Gamma(p_{t+1})>\delta>0,
$$
and if $\Gamma$ is a broken tree then $\phi_\Gamma$ satisfies the join equation \eqref{e:join}.
\end{lma}

\begin{pf}
Possible areas of trees of $L$ with one positive puncture constitutes a finite set of numbers
$$
0<\alpha_1<\alpha_2\dots<\alpha_m.
$$
Let ${\overline \M}_j$ denote the compactification of the moduli space of trees of area $\le\alpha_j$. Pick a stretching function $\Gamma\mapsto\phi_\Gamma$ for $\Gamma\in\M_1$. Assume inductively that we have defined a function ${\overline \M_j}\to C^k(S^1,\R)$ in such a way that there exists an $\epsilon>0$ with the following properties. In an $\epsilon$-neighborhood of each broken tree the function satisfies the join equation. For $\Gamma$ outside a $2\epsilon$-neighborhood of the boundary of ${\overline \M}_j$ the functions $\phi_\Gamma$ are stretching. In the region between these neighborhoods, the inequality,
$$
\phi_\Gamma(p_t)-\phi_\Gamma(p_{t+1})
>\frac12\left(1-\frac{L}{(K+1)}\right)\frac{1}{M}h(l(p_t,p_{t+})),
$$
where the notation is as in Lemma \ref{l:jumpbroken}, holds. We want to extend the function to ${\overline \M}_{j+1}$. Since the boundary of ${\overline \M}_{j+1}$ consists of broken trees of area strictly smaller than $\alpha_{j+1}$ we define the map $\pa{\overline \M}_{j+1}\to C^k(S^1,\R)$ by imposing the join equation in a small neighborhood of the boundary of ${\overline \M}_{j+1}$. Using Lemma \ref{l:jumpbroken} it is not hard to see that if $\epsilon>0$ is sufficiently small then we can extend this family to all of ${\overline \M}_{j+1}$ respecting the above conditions on $\epsilon$-neighborhood as well as on the $2\epsilon$-neighborhood and its complement. We obtain a function on ${\overline \M}$ with properties as desired.
\end{pf}

We view a function ${\overline \M}\to C^{k}(S^{1},\R)$ with properties as in Lemma \ref{l:shiftfct} as a function into $C^{k}(S^{1}, T^{\ast} I)$ by identifying $\R$ with the fiber in $T^{\ast}I$. Fix such a function and use it as a perturbation function for flow trees, see Subsection \ref{s:pert}. Let $\Delta_{\rm a}$ denote the corresponding differential on $\A(L\times I)$.

\begin{rmk}\label{r:smallvar}
By scaling, the total variation of any function $\phi_\Gamma$, $\Gamma\in{\overline \M}$ can be assumed arbitrarily small.
\end{rmk}

\begin{lma}\label{l:onelevel}
Any rigid perturbed flow tree has only one level. The flow tree corresponding to the perturbed slice tree in that level is itself rigid and exactly one of its negative punctures lies at the $0$--level, or the tree is entirely contained in $\{x_0=\pm 1\}$.
\end{lma}

\begin{pf}
This follows from the choice of abstract perturbation. Since it orders the negative punctures, at most one at a time can lie at the $0$--level. The lemma is then an easy consequence of the dimension formula.
\end{pf}

Let $\pa\colon\A(L)\to\A(L)$ denote the contact homology differential and let the notation for generators be $\hat c$ for the maximum in $I_c$ and $c[\pm 1]$ for the minima at $\pm 1\in I_c$, for any Reeb chord $c$ of $L$.

\begin{cor}\label{c:abstDGA}
The differential $\Delta_{\rm a}\colon \A(L\times I)\to\A(L\times I)$ satisfies the following
\begin{align}\label{e:pm1}
\Delta_{\rm a} c[\pm1] &=\pa c[\pm1],\\\label{e:hat}
\Delta_{\rm a} \hat c&=c[+1]+c[-1]+\Gamma(\pa c),
\end{align}
where $\Gamma(1)=0$ and 
\begin{align*}
\Gamma(a_1 a_2\dots a_m)&=\hat a_1a_2[-1]\dots a_m[-1]\\
                        &+a_1[+1]\hat a_2a_3[-1]\dots a_m[-1]\\
                        &\quad\vdots\\
                        &+a_1[+1]\dots a_{m-1}[+1]\hat a_m.
\end{align*}
\end{cor}
\begin{pf}
It is easy to construct the rigid perturbed trees mentioned in Lemma \ref{l:onelevel}, which give the second term in \eqref{e:hat}, as well as the flow lines of $\beta_c$ which give the first term in \eqref{e:hat}. Lemma \ref{l:onelevel} shows that there are no other rigid perturbed trees.
\end{pf}

Let $\Delta\colon\A(L\times I)\to\A(L\times I)$ denote the differential on the contact homology algebra of the trace of the constant isotopy which arises from a geometric deformation.
\begin{cor}
There exists a tame isomorphism between the algebras $(\A(L\times I),\Delta)$ and $(\A(L\times I),\Delta_{\rm a})$.
\end{cor}

\begin{pf}
Remark \ref{r:geomasabstr} shows that $\Delta$ can be expressed as a differential $\Delta_{\rm g}$ induced by a perturbation function and by collections of Morse functions on the Reeb chord manifolds. Corollary \ref{c:tameiso} then shows that $(\A(L\times I),\Delta_{\rm a})$ and $(\A(L\times I),\Delta_{\rm g})$ are tame isomorphic.
\end{pf}

\section{Decomposing isotopies and differentials}\label{S:move}
In this section we first discuss how to subdivide a Legendrian isotopy into pieces in such a way that each piece corresponds either to a sufficiently good approximation of a constant isotopy or to a move isotopy which is constant except in a small region where it has one of several standard forms. The computations of Section \ref{S:abstrprt} allow us to describe the DGA of the trace of the almost constant isotopy, see Subsection \ref{s:nearconst}. Slight extensions of these computations allow us to describe the DGAs of traces of move isotopies as well, see Subsections \ref{s:frontmove} and  \ref{s:lagmove}.

\subsection{The trace of an almost constant isotopy}\label{s:nearconst}
Let $L\subset J^{1}(\R)$ be a generic Legendrian submanifold. Let $L_t$, $-1\le t\le 1$ be an isotopy of $L=L_0$ and let $\Phi(L\times I)$ be the trace of $L_t$ with standard ends. Let $\Psi\colon J^{1}(\R\times I)\to J^{1}(\R\times I)$ be the perturbation which makes the trace $L\times I$ of the constant isotopy generic with respect to flow trees.
\begin{lma}\label{l:almostconst}
There exists $\epsilon=\epsilon(L)>0$ such that if the isotopy $L_t$ is contained in a the $C^{2}$ $\epsilon$-neighborhood of the constant isotopy then the contact homology algebra of $\Psi(\Phi(L\times I))$ and that of $\Psi(L\times I)$ are canonically isomorphic.
\end{lma}

\begin{pf}
Assume not. Then there exists a sequence of isotopy traces $\Phi_j\colon L\times I\to J^{1}(\R\times I)$, $j=1,2,\dots$ such that the space of rigid flow trees of $\Psi(\Phi_j(L\times I))$ and that of $\Psi(L\times I)$ are non-isomorphic and such that $\Phi_j\to\id$ as $j\to\infty$. By the compactness properties of the space of flow trees, this contradicts the flow trees of $\Psi(L\times I)$ being transversely cut out .
\end{pf}

Let $L_t$, $0\le t\le 1$ be a Legendrian isotopy. It is well known that we may deform this isotopy into an isotopy with the property that $L_t$ is a generic Legendrian submanifold for $t\ne t_j$ where $0<t_1<\dots< t_m<1$ is a finite collection of instances and such that around every $t_j$ the isotopy is constant outside a small disk and inside this disk the Legendrian undergoes one of the following standard moves. We call such an isotopy an isotopy with standard moves.

\begin{itemize}
\item[({\bf F1})] A first Redemeister front move.
\item[({\bf F2})] A second Redemeister front move.
\item[({\bf F3})] A third Redemeister front move.
\item[({\bf L1})] A Lagrangian triple point move.
\item[({\bf L2})] A pair of Reeb chords disappear.
\item[({\bf L3})] A pair of Reeb chords appear.
\end{itemize}

Given an isotopy with standard moves we may subdivide the intervals between the standard moves in sufficiently short isotopies through generic Legendrian submanifolds. This gives a representative of the trace of the isotopy which is a join of approximately constant isotopies and isotopies with standard moves. Let $L$ be the middle stage of an almost constant isotopy and let $\{\hat c, c[-1], c[+1]\}$ denote the Reeb chords of the trace of the almost constant isotopy, where $c$ ranges over the Reeb chords of $L$, where $\hat c$ corresponds to the maximum in $I_c$ and $c[\pm 1]$ to the minima at the endpoints of $I_c$. If $w\in\A(L)$ is an element then write $w[\pm 1]$ for the corresponding element in $\A(L\times I)$ where every generator $c$ has been replaced by $c[\pm 1]$. 

\begin{cor}\label{c:acdiff}
The DGA differential $\Delta$ of the trace of an almost constant isotopy is the following
\begin{align}\label{e:iddiff}
&\Delta c[\pm 1]=(\pa c)[\pm 1],\\\notag
&\Delta\hat c = c[-1]+(\id c)[+1] + \Gamma_{\id}(\pa c),
\end{align}
for every Reeb chord $c$.
\end{cor}

\begin{pf}
This follows from Lemma \ref{l:diffconst} in combination with Lemma \ref{l:almostconst}.
\end{pf}

Concatenation of the almost constant isotopies of an isotopy without moves then gives an expression for the differential of the trace of the total isotopy.

\subsection{Front moves}\label{s:frontmove}
In this subsection we study moves {\rm ({\bf F1})}--{\rm ({\bf F3})}. In fact, none of these moves change the differential from that of the trace of the constant isotopy.

\begin{lma}\label{l:F1}
The DGA of the trace of an {\rm ({\bf F1})}-move is identical to the DGA of the trace of an almost constant isotopy, see \eqref{e:iddiff}.
\end{lma}

\begin{figure}
\centering
\includegraphics[width=.8\linewidth]{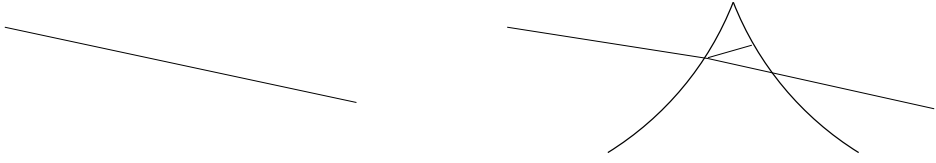}
\caption{Flow line of the trace of the constant isotopy (left) and corresponding flow tree of the {\rm ({\bf F1})}-isotopy (right).}
\label{fig:F1}
\end{figure}

\begin{pf}
We parametrize the trace by keeping the isotopy constant on $x_0\in [-1,\frac12]$, then performing the move, and keeping it again constant close to $x_0=1$. It is straightforward to check that there is a natural 1--1 correspondence between flow trees of the trace of the constant isotopy and flow trees of the {\rm ({\bf F1})}-isotopy. The only difference between the two kinds of trees appears as they pass the bifurcation region where trees in one of the situations are obtained from trees in the other situation by adding a $Y_1$-vertex and an end, see Figure \ref{fig:F1}.
\end{pf}

\begin{lma}\label{l:F2}
The DGA of the trace of an {\rm ({\bf F2})}-move is identical to the DGA of the trace of an almost constant isotopy, see \eqref{e:iddiff}.
\end{lma}

\begin{pf}
As in the proof of Lemma \ref{l:F1} it is straightforward to find a 1--1 correspondence between rigid flow trees of the trace of the constant isotopy and the trace of the isotopy of the move. To see this, use the cotangent lift of the trees and the fact that the Lagrangian projection of the trace of an {({\bf F2})}-isotopy is qualitatively indistinguishable from that of the constant isotopy.
\end{pf}

\begin{lma}\label{l:F3}
The DGA of the trace of an {\rm ({\bf F3})}-move is identical to the DGA of the trace of an almost constant isotopy, see \eqref{e:iddiff}.
\end{lma}

\begin{pf}
The same argument as in the proof of Lemma \ref{l:F2} applies.
\end{pf}

\subsection{Lagrangian moves}\label{s:lagmove}
In this subsection we study the moves {\rm ({\bf L1})}--{\rm ({\bf L3})}. To this end we will append the move isotopy to one end of a constant isotopy. As in Subsection \ref{s:pert} our computation uses abstract perturbations. That is, we use perturbed flow trees in our computations. The argument which shows that this computation gives a DGA which is tame isomorphic to the DGA which arises from a geometric perturbation is almost identical to the proof of Lemma \ref{l:diffconst} and Corollary \ref{c:tameiso} and will not be repeated. 

Let $L_{-1}$ and $L_{+1}$ be Legendrian submanifolds of $J^{1}(\R)$ such that $L_{+1}$ is obtained from $L_{-1}$ via an {\rm ({\bf Lj})}-move. Then we write $\phi_{{\rm ({\bf Lj})}}\colon \A(L_{-1})\to\A(L_{+1})$ for the induced homomorphism, ${\bf j}={\bf 1},{\bf 2},{\bf 3}$, see \cite{K}. Furthermore we write $\pa_\pm$ for the differential on $\A(L_{\pm 1})$ and if $\phi\colon\A(L_{-1})\to\A(L_{+1})$ is a homomorphism, we use the notion $\Gamma_\phi\colon\A(L_{-1})\to\A(L\times I)$ as in Theorem \ref{t:main}.

We first consider the simpler cases of {\rm ({\bf L1})} and {\rm ({\bf L2})} where the abstract perturbations and the move region can be taken disjoint. We let the move happen inside a box of the form $\frac12\le x_0\le \frac12+\delta$, $|x_1-a|\le \epsilon$. Inside this box we can draw the flow explicitly and from that information compute the differential.

\begin{lma}\label{l:L1}
The differential $\Delta$ of the trace of a move isotopy of type {\rm ({\bf L1})} satisfies
\begin{equation}\label{e:L1}
\Delta \hat c = c[-1]+\phi_{\bf L1}(c)[+1]+\Gamma_{\phi_{\bf L1}}(\pa_- c)
\end{equation}
\end{lma}

\begin{pf}
We chose abstract perturbations of the same form as in Subsection \ref{s:pert} and we use a cut-off function supported above the move region. As there, we find that as $\eta\to 0$ outside the move region any sequence of trees converges to a generalized tree. In particular, any rigid flow tree has all its limit slice trees near $x_0=0$. The difference arises as the descending flow lines near the Reeb chord manifolds enter the box of the move. There are two different {\rm ({\bf L1})}-moves see Figures \ref{fig:L1i1} and \ref{fig:L1ii1}. The corresponding flow pictures are shown in Figures \ref{fig:L1i2} and \ref{fig:L1ii2}, respectively. It is clear from these pictures that there is no effect on the differential in the first case and that in the second case we get exactly $c[-1]+\phi(c)[+1] + \Gamma_\phi(\pa c)$ in the right hand side of \eqref{e:L1} where $\phi(a)=a+bc$ and $\phi$ is the identity on other generators.
\begin{figure}
\centering
\includegraphics[width=.8\linewidth]{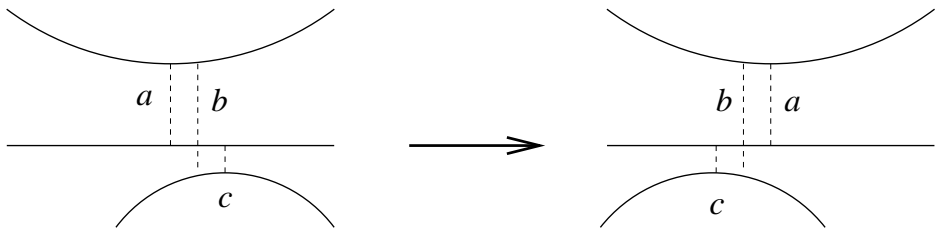}
\caption{The first {\rm ({\bf L1})}-isotopy.}
\label{fig:L1i1}
\end{figure}

\begin{figure}
\centering
\includegraphics[width=.8\linewidth]{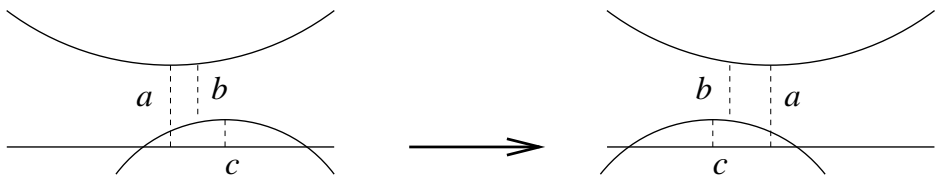}
\caption{The second {\rm ({\bf L1})}-isotopy.}
\label{fig:L1ii1}
\end{figure}

\begin{figure}
\centering
\includegraphics[width=.4\linewidth]{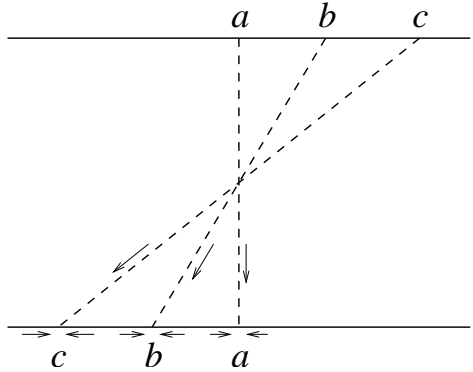}
\caption{The trace of the first {\rm ({\bf L1})}-isotopy. Every flow line belonging to a rigid tree passes right through the isotopy box.}
\label{fig:L1i2}
\end{figure}

\begin{figure}
\centering
\includegraphics[width=.4\linewidth]{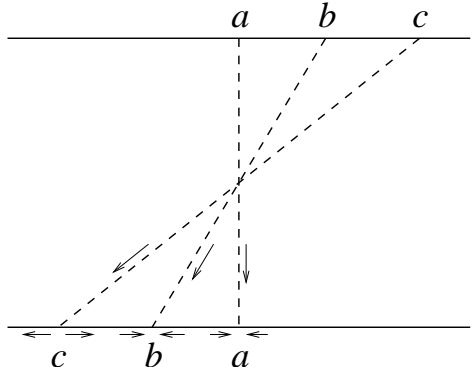}
\caption{The trace of the second {\rm ({\bf L1})}-isotopy. Any flow line on its way toward $a$ can either split over a $Y_0$-vertex into a rigid configuration of a flow line toward $b$ and a flow line toward $c$ or pass right through the isotopy box.}
\label{fig:L1ii2}
\end{figure}
\end{pf}

\begin{lma}\label{l:L2}
The differential $\Delta$ of the trace of a move isotopy of type {\rm ({\bf L2})} satisfies
$$
\Delta \hat c = c[-1]+\phi(c)[+1]+\Gamma_{\phi_{\bf L2}}(\pa_- c)
$$
\end{lma}

\begin{pf}
The proof is very similar to the proof of Lemma \ref{l:L1}. We use abstract perturbations supported outside the move box. The move is depicted in Figure \ref{fig:L2f} and the corresponding flow box in Figure \ref{fig:L2}.
\begin{figure}
\centering
\includegraphics[width=.8\linewidth]{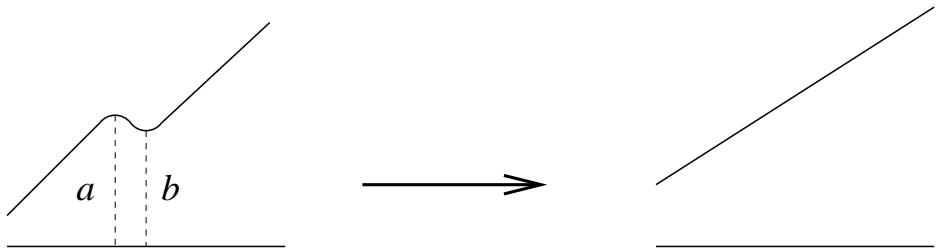}
\caption{An {\rm ({\bf L2})}-isotopy.}
\label{fig:L2f}
\end{figure}
\begin{figure}
\centering
\includegraphics[width=.4\linewidth]{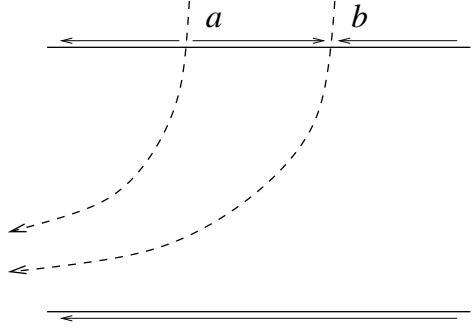}
\caption{The trace of the {\rm ({\bf L2})}-isotopy. Any flow line on its way toward $a$ goes on toward a word contributing to $v$. The corresponding tree lives in a $1$-parameter family since nearby flow lines (which would formerly be attracted by $b$) also goes on toward the same word in $v$. A flow line of a rigid tree on its way toward $b$ is determined uniquely already in the region above. It goes on to a word contributing to $v$.}
\label{fig:L2}
\end{figure}
It follows that
\[
\Delta\hat c = c[-1] + \phi_{{\bf L2}}(c)[+1]+\Gamma_{\phi_{{\bf L2}}}(\pa_- c)
\]
where $\phi$ takes $a$ to $0$, $b$ to $\pa a+b=v$, and is the identity on all other generators.
\end{pf}

\begin{lma}\label{l:L3}
The differential $\Delta$ of the trace of a move isotopy of type {\rm ({\bf L3})} satisfies
$$
\Delta\hat c = c[-1]+\phi_{\bf L3}(c)[+1]+\Gamma_{\phi_{{\bf L3}}}(\pa c).
$$
\end{lma}

\begin{pf}
To prove this result we use an abstract perturbation of the kind discussed in Section \ref{S:abstrprt}, which orders the punctures according to their natural boundary ordering. This time we let the perturbation of the horizontal trees be cut-off only after the flow box of the move. This flow box is depicted in Figure \ref{fig:L3}.
\begin{figure}
\centering
\includegraphics[width=.4\linewidth]{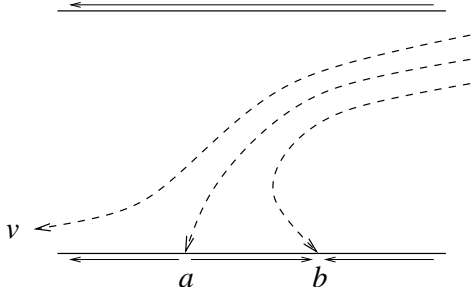}
\caption{The flow of an {({\bf L3})}-isotopy. Consider a generalized perturbed flow tree near the $[-1]$-chords moving upwards. When one of its $b$-chords hits the flow line going to $a$ in the flow box, a $1$-parameter family of perturbed flow trees becomes rigid. When the puncture moves further up the flow goes instead to a word contributing to $v$ and stays there until the next $b$-puncture hits the flow line ending at $a$.}
\label{fig:L3}
\end{figure}
Our choice of perturbation guarantees that the flow lines of a generalized perturbed tree hits the $a$-flow line in the order they appear on the boundary. The formula for the differential follows.
\end{pf} 
\section{Algebraic treatment of concatenation and of homotopy of isotopies}\label{S:algebra}
In this section we first show that the DGAs of traces, after destabilization and tame isomorphisms have certain naturality properties with respect to concatenation. Together with the results of Section \ref{S:move} this leads to a proof of Theorem \ref{t:main}. Second we show that the contact homology of the trace of an isotopy depends only on the chain homotopy type of the contact homology morphism induced by the isotopy. Geometrically, such a chain homotopy corresponds to a homotopy of isotopies.

\subsection{Concatenation and the proof of Theorem \ref{t:main}}\label{s:algconc}
\begin{figure}
\labellist
%\small\hair 2pt
\pinlabel $\A_-$ at 6 152
\pinlabel $\A_+$ at 6 10
\pinlabel $\varphi$ at -10 80
\pinlabel $\iota_-$ at 65 165
\pinlabel $\iota_+$ at 65 20
\pinlabel $\Gamma_\varphi$ at 60 105
\pinlabel $\CC_\varphi$ at 280 80
\endlabellist
\centering
\includegraphics[width=.4\linewidth]{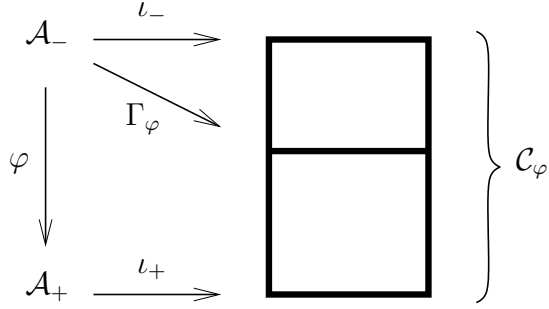}
\caption{The mapping cone of a DGA morphism}
\label{fig:cone}
\end{figure}

Let $(\A_\pm,\pa_\pm)$ be a DGA with a finite set of generators $g(\A_\pm)$ and let $\varphi\colon\A_-\to\A_+$ be a DGA morphism. Let the non-commutative algebra $\CC_\varphi$ be generated by the following:
\begin{itemize}
\item a generator $a_\pm$ for each $a_\pm\in g(\A_\pm)$, with grading $|a_\pm|$ as in $\A_\pm$,
\item a generator $\hat a$ for each $a_-\in g(\A_-)$, with grading $|\hat a|=|a_-|+1$.
\end{itemize}
There are obvious embeddings (algebra monomorphisms) $\iota_\pm\colon\A_\pm\to\CC_\varphi$ defined (on generators) by $\iota_\pm(a_\pm)=a_\pm$, see Figure \ref{fig:cone}.

We define a map $\Delta_\varphi\colon \CC_\varphi\to\CC_\varphi$ as follows on generators
\begin{align*}
\Delta_\varphi a_\pm &= \iota_\pm(\pa_\pm a_\pm),\\
\Delta_\varphi\hat a &= a_- +\iota_+(\varphi(a_-))+\Gamma_\varphi(\pa_- a_-)
\end{align*}
and extend it to all of $\CC_\varphi$ by the Leibniz rule. Here $\Gamma_\varphi\colon\A_-\to\CC_\varphi$ is defined as follows. On generators $a_-$ of $\A_-$,
\[
\Gamma_\varphi(a_-)=\hat a
\]
and for products $uv\in\mathcal A_-$,
\[
\Gamma_\varphi(uv)=\iota_-(u)\Gamma_\varphi(v)+\Gamma_\varphi(u)\iota_+(\varphi(v)).
\]
The $(\iota_-,\iota_+\circ\varphi)$--derivation $\Gamma_\varphi$ increases grading by $1$ and consequently $\Delta_\varphi$ decreases grading by $1$.

Since $\pa_\pm$ also satisfies the Leibniz rule, it follows that $\Delta_\varphi \iota_\pm(u_\pm)=\iota_\pm(\pa_\pm u_\pm)$ for all $u_\pm\in\A_\pm$. In other words, $\iota_\pm$ will be DGA morphisms once we show that $(\CC_\varphi,\Delta_\varphi)$ is a DGA. To establish that we must show that $\Delta_\varphi^2=0$ and it is enough to prove this for generators of $\CC_\varphi$. For generators $a_\pm$ this is trivial and for generators $\hat a$, it follows from the following lemma.

\begin{lma}
For all $u\in\A_-$,
\begin{equation}\label{e:Delta^2=0}
D(u)=\iota_-(u)+\iota_+(\varphi(u))+\Delta_\varphi\Gamma_\varphi(u)+\Gamma_\varphi(\partial_- u)=0.
\end{equation}
In particular $\Delta_\varphi^2\hat a=0$ for any generator $\hat a$.
\end{lma}

\begin{pf}
If $a_-$ is a generator of $\A_-$, then
\begin{align*}
D(a_-)&=\iota_-(a_-)+\iota_+(\varphi(a_-))+\Delta_\varphi\Gamma_\varphi(a_-)+\Gamma_\varphi(\partial_- a_-)\\
&=a_-+\iota_+(\varphi(a_-))+\Delta_\varphi\hat a+\Gamma_\varphi(\partial_- a_-)=0
\end{align*}
by definition of $\Delta_\varphi\hat a$. If $D(u)=D(v)=0$, then
\begin{align*}
D(uv)&=\iota_-(uv)+\iota_+(\varphi(uv))+\Delta_\varphi\Gamma_\varphi(uv)+\Gamma_\varphi(\partial(uv))\\
&=\iota_-(u)\iota_-(v)+\iota_+(\varphi(u))\iota_+(\varphi(v))
+\Delta_\varphi\left(\iota_-(u)\Gamma_\varphi(v)+\Gamma_\varphi(u)\iota_+(\varphi(v))\right)\\
&+\Gamma_\varphi(\partial_- u\cdot v+u\cdot\partial_- v)\\
&=\iota_-(u)\iota_-(v)+\iota_+(\varphi(u))\iota_+(\varphi(v))\\
&+\Delta_\varphi(\iota_-(u))\cdot\Gamma_\varphi(v)+
\iota_-(u)\cdot\Delta_\varphi\Gamma_\varphi(v)+\Delta_\varphi\Gamma_\varphi(u)\cdot \iota_+(\varphi(v))
+\Gamma_\varphi(u)\cdot\Delta_\varphi \iota_+(\varphi(v))\\
&+\iota_-(\partial_- u)\Gamma_\varphi(v)+\Gamma_\varphi(\partial_- u)\iota_+(\varphi(v))+\iota_-(u)\Gamma_\varphi(\partial v)+\Gamma_\varphi(u)\iota_+(\varphi(\partial v))\\
&=\left(\iota_+(\varphi(u))+\Delta_\varphi\Gamma_\varphi(u)+\Gamma_\varphi(\partial_- u)\right)j(\varphi(v))+\iota_-(u)\left(\iota_-(v)+\Delta_\varphi\Gamma_\varphi(v)+\Gamma_\varphi(\partial_- v)\right)\\
&=\iota_-(u)\iota_+(\varphi(v))+\iota_-(u)\iota_+(\varphi(v))=0.\\
\end{align*}
For the second statement, note that $\Delta^{2}_\varphi\hat a= D(\pa_- a)$.
\end{pf}

We call the DGA $(\CC_\varphi,\Delta_\varphi)$ the {\em mapping cone} of the chain map $\varphi\colon\A_-\to\A_+$.

The proof of the main theorem of this section relies on the following result which is often used in the subject of contact homology. For a proof see \cite{Ch}.

\begin{thm}\label{thm:stab}
Let $\A=T(q_1,q_2,\ldots,q_m,a,b)$ and $\A'=T(q_1,q_2,\ldots,q_m)$ be DGAs freely generated by the indicated generators with differentials $\partial$ and $\partial'$, respectively, which decrease grading by one and which satisfy the Leibniz rule. Assume that $\A$ and $\A'$ come equipped with height filtrations, that is that there exists an ordering
of the generators $c\in\{q_1,\dots,q_m,a,b\}$ of $\A$ with the following property. For each $c$, the expression $\partial c$ ($\partial'c$) is a polynomial of other generators of $\A$ (of $\A'$) in which all generators is of smaller height than $c$. Assume that $\partial a=b+v$. Define the projection
\[
\tau\colon\A\to\A'
\]
by $\tau(q_j)=q_j$, $1\le j\le m$, $\tau(a)=0$, and $\tau(b)=v$.
Extend $\tau$ to a grading preserving algebra homomorphism and suppose it is a chain map.
Then, the DGAs $(\A,\partial)$ and $(\A',\partial')$ are stable tame isomorphic. In particular they have isomorphic homologies.
\end{thm}

\begin{figure}
\labellist
%\small\hair 2pt
\pinlabel $\mathcal A_1$ at 10 260
\pinlabel $\mathcal A_2$ at 10 135
\pinlabel $\mathcal A_3$ at 10 10
\pinlabel $\alpha$ at -10 195
\pinlabel $\beta$ at -10 70
\pinlabel $\Gamma_\alpha$ at 50 220
\pinlabel $\Gamma_\beta$ at 50 90
\pinlabel $k$ at 55 25
\pinlabel $j$ at 55 150
\pinlabel $i$ at 55 275
\pinlabel $\tau_1$ at 275 150
%\pinlabel $\psi_1$ at 275 95
\pinlabel $\tau_2$ at 490 150
%\pinlabel $\psi_2$ at 490 95
\pinlabel $\tau_k$ at 635 150
%\pinlabel $\psi_k$ at 635 95
\pinlabel $\mathcal B=\mathcal B_0$ at 170 280
\pinlabel $\mathcal B_1$ at 385 280
\pinlabel $\mathcal B_k=\mathcal C_{\beta\alpha}$ at 745 280
\endlabellist
\centering
\includegraphics[width=.8\linewidth]{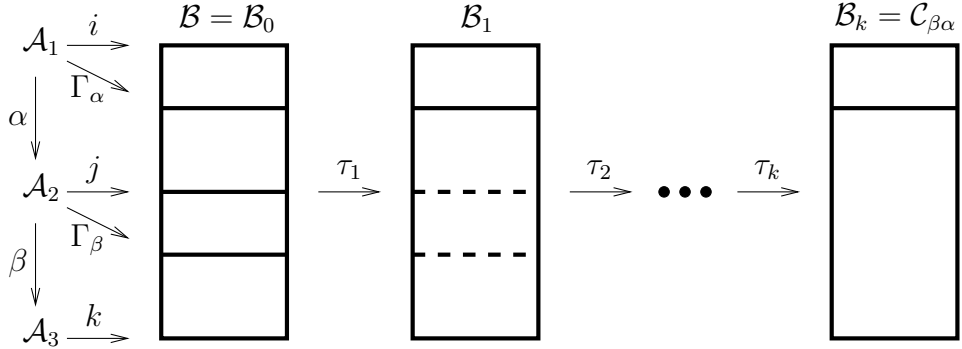}
\caption{Concatenating mapping cones}
\label{fig:comp}
\end{figure}

Let now $(\A_1,\partial_1)$, $(\A_2,\partial_2)$, and $(\A_3,\partial_3)$ be filtered DGAs with $\alpha\colon\A_1\to\A_2$ and $\beta\colon\A_2\to\A_3$ DGA morphisms between them (see Figure \ref{fig:comp}). Let $(\CC_\alpha,\Delta_\alpha)$ be the mapping cone of $\alpha$ and $(\CC_\beta,\Delta_\beta)$ be the mapping cone of $\beta$. Define the algebra $\mathcal B$ by taking the disjoint union of $\CC_\alpha$ and $\CC_\beta$ and identifying for each generator $b$ of $\mathcal A_2$, the corresponding generator $b\in\CC_\alpha$ with $b\in\CC_\beta$. This is the algebraic version of concatenation and $\Delta_\alpha$ and $\Delta_\beta$ define a differential $\Delta$ on $\mathcal B$.
Note that the inclusions of $\A_j\to\mathcal B$, $j=1,2,3$ are DGA morphisms.

\begin{lma}\label{l:algcomp}
The DGA $(\mathcal B,\Delta)$ is stable tame isomorphic to the mapping cone $(\CC_{\beta\alpha},\Delta_{\beta\alpha})$ of the chain map $\beta\circ\alpha\colon\A_1\to\A_3$.
\end{lma}

%%%%%%%%%%%%%
%% remove the psis! %%
%%%%%%%%%%%%%

\begin{pf}
We remove pairs of generators $(b,\hat b)$ for $b\in\mathcal A_2$ inductively. Let the generators of $\mathcal A_2$ be arranged by height from $b_1$ (the longest) to $b_k$ (the shortest). We eliminate $b_1$ first, then the second longest generator $b_2$, and so on until we reach $b_k$. Starting from $\mathcal B=\mathcal B_0$, we obtain the DGAs $\mathcal B_i$ ($i=0,1,\ldots,k$) that are generated by $a$ and $\hat a$ for each generator $a\in\A_1$, $c$ for each generator $c\in\A_3$, and by $b_{i+1},\hat b_{i+1}\ldots,b_k,\hat b_k$. We may assume that in each $\mathcal B_i$, the heights of the $\hat a$ generators are all higher than those of the $b$ and $\hat b$ generators.

Starting from $\Delta_0=\Delta$, we define differentials $\Delta_i$ on $\mathcal B_i$ recursively. If $\Delta_{i-1}\hat b_i=b_i+v_i$ then we define the projection $\tau_i\colon\mathcal B_{i-1}\to\mathcal B_i$ as in Theorem \ref{thm:stab}. (We remark that $\Delta_{i-1}\hat b_i$ will always have that form, in fact with $v_i=\beta(b_i)+\Gamma_\beta(\partial_2b_i)$, since we eliminate in the order of decreasing height; the interesting part of this definition is $\Delta_i\hat a$ for $a\in\mathcal A_1$.) For any generator $c$ of $\mathcal B_i$, we define $\Delta_ic=\tau_i(\Delta_{i-1}c)$. Then $\tau_i$ is a DGA morphism and Theorem \ref{thm:stab} implies that all $\mathcal B_i$ are stable tame isomorphic.

To finish the proof we use the following: for all $b\in\A_2\subset \mathcal B$,
\begin{equation}\label{eq:tech}
\tau_k\circ\tau_{k-1}\circ\cdots\circ\tau_1(b)=\beta(b).
\end{equation}
To see this, note that
\[
\tau_k\circ\tau_{k-1}\circ\cdots\circ\tau_1(b_i)=\tau_k\circ\cdots\circ\tau_i(b_i)=\tau_k\circ\cdots\circ\tau_{i+1}(\beta(b_i)+\Gamma_\beta(\partial_2 b_i)).
\]
Here, each monomial of $\Gamma_\beta(\partial_2b_i)$ contains a $\hat b_j$-factor for some $i+1\le j\le k$, thus one of the remaining projections will annihilate it, and $\beta(b_i)$ is in $\A_3$ and is therefore fixed by all projections.

We claim that $(\mathcal B_k,\Delta_k)$ is isomorphic to the mapping cone $(\mathcal C_{\beta\alpha},\Delta_{\beta\alpha})$. Their sets of generators agree and the differentials for generators $a\in\mathcal A_1$ and $c\in\mathcal A_3$ are obviously the same. For any generator $a\in\mathcal A_1$,
\begin{eqnarray*}
\Delta_k\hat a&=&\Delta_k\tau_k\circ\tau_{k-1}\circ\cdots\circ\tau_1(\hat a)\\
&=&\tau_k\circ\tau_{k-1}\circ\cdots\circ\tau_1(\Delta_0\hat a)\\
&=&\tau_k\circ\tau_{k-1}\circ\cdots\circ\tau_1(a+\alpha(a)+\Gamma_\alpha(\partial_1a))\\
&=&a+\tau_k\circ\tau_{k-1}\circ\cdots\circ\tau_1(\alpha(a))+\tau_k\circ\tau_{k-1}\circ\cdots\circ\tau_1(\Gamma_\alpha(\partial_1 a))\\
&=&a+\beta(\alpha(a))+\Gamma_{\beta\circ\alpha}(\partial_1a)).
\end{eqnarray*}
To see that the last equality holds, note that each summand of $\Gamma_\alpha(\partial_1 a)$ is a three-fold product of an element of $\mathcal A_1$ (possibly $1$), an $\hat a$--generator (these are all fixed by projections) and the $\alpha$--image of an element of $\mathcal A_1$, note that projections are algebra homomorphisms and use \eqref{eq:tech}.
\end{pf}

\begin{pf}[Proof of Theorem \ref{t:main}]
The theorem follows from Lemma \ref{l:almostconst}, Corollary \ref{c:acdiff}, Lemmas \ref{l:F1} -- \ref{l:L3}, and Lemma \ref{l:algcomp}.
\end{pf}

\subsection{Chain homotopy}\label{s:chiso}
Consider two DGAs $(\A_+,\pa_+)$ and $(\A_-,\pa_-)$. Let $\psi,\phi\colon\A_+\to\A_-$ be chain maps and let $K\colon \A_+\to\A_-$ be a chain homotopy between them. That is,
\[
\phi+\psi=K\circ\pa_+ + \pa_-\circ K.
\]
Consider the mapping cone DGAs $(\CC_\phi,\Delta_\phi)$ and $(\CC_\psi,\Delta_\psi)$. Define the map $\Gamma_K\colon \A_+\to\CC_\psi$ to equal zero on constants and on linear monomials and by the following expression for monomials $b_1\dots b_r$ of length $r\ge 2$
\[
\Gamma_K(b_1\dots b_r) = \hat b_1 K(b_2\dots b_r) + b_1\hat b_2 K(b_3\dots b_r) + \dots +
b_1\dots b_{r-2}\hat b_{r-1}K(b_r).
\]

\begin{lma}
The algebra map $F\colon (\CC_\phi,\Delta_\phi)\to (\CC_\psi,\Delta_\psi)$ defined on generators as follows
\begin{align*}
&F(c)=c, \qquad c\in\A_+,\\
&F(v)=v, \qquad v\in\A_-,\\
&F(\hat c)=\hat c + K(c) + \Gamma_K(\pa_+ c), \qquad c\in\hat\A_+,
\end{align*}
is a (tame) chain isomorphism. That is,
\[
F\circ\Delta_\phi=\Delta_\psi\circ F.
\]
\end{lma}

\begin{pf}
For simpler notation we write $\CC_\phi=\A=\CC_\psi$ and consider this as one algebra with two differentials and in order to facilitate computations we introduce the following notation. Let $H\colon\A_+\to\A$ and $\theta\colon \A_+\to A$ be maps. Then define the map $\Omega_\theta^{H}\colon\A_+\to\A$ as follows on monomials
\[
\Omega_\theta^{H}(b_1\dots b_r)= H(b_1)\theta(b_2\dots b_r) + b_1 H(b_2)\theta(b_3\dots b_r) + \dots + b_1\dots b_{r-1}H(b_r).
\]
We compute
\begin{align}\label{e:comp1}
F(\Delta_\phi(\hat c)) &=
F(c + \phi(c) + \Gamma_\phi(\pa_+ c))\\\notag
&=c + \phi(c) + \Gamma_\phi(\pa_+ c) + \Omega_\phi^{K}(\pa_+ c) + \Omega_\phi^{\Gamma_K\circ\pa_+}(\pa_+ c),
\end{align}
and
\begin{align}\label{e:comp2}
\Delta_\psi(F(\hat c)) &=
\Delta_\psi(\hat c + K(c) + \Gamma_K(\pa_+ c))\\\notag
&=c + \psi(c) + \Gamma_\psi(\pa_+ c) + \pa_-(K(c)) + \Delta_\psi(\Gamma_K(\pa_+ c)).
\end{align}
We note that the monomials in \eqref{e:comp1} and \eqref{e:comp2} are of two types: monomials which are constant in $\hat\A_+$-generators and monomials which are linear in $\hat\A_+$-generators. We first show that the monomials of the former kind cancels between the two equations.

The contribution from \eqref{e:comp1} to monomials of the first kind is
\[
c+\phi(c)+\Omega_\phi^{K}(\pa_+ c)
\]
and the contribution from \eqref{e:comp2} is
\[
c+\psi(c)+\pa_-(K(c))+\Delta_\psi(\Gamma_K(\pa_+ c))_0,
\]
where $w_0$ denotes the term in an element which is constant in the $\hat\A_+$-generators.
Thus, if we show that
\begin{equation}\label{e:monomial1}
\Omega_\phi^{K}(\pa_+ c)+\Delta_\psi(\Gamma_K(\pa_+ c))_0=K(\pa_+ c)
\end{equation}
then it follows that monomials of the first type cancels. Now, if $b_1\dots b_r$ is a monomial in $\pa_+ c$ then the contribution from this monomial to the left hand side of \eqref{e:monomial1} which is non-constant in the $\A_+$-generators vanishes for the following reason. Terms arising from $\Omega_\phi^{K}$ have the form
\[
b_1\dots b_s K(b_{s+1})\phi(b_{s+2})\dots\phi(b_r)
\]
and are canceled by terms in $\Delta_\psi(b_1\dots \hat b_s K(b_{s+1})\phi(b_{s+2})\dots\phi(b_{r}))$ corresponding to the $b_s$-term in $\Delta_\psi(\hat b_s)$. Remaining terms in $\Delta_\psi(b_1\dots \hat b_s K(b_{s+1}\dots b_{r}))$ corresponding to the $b_s$-term in $\Delta_\psi(\hat b_s)$ cancels with terms in $\Delta_\psi(b_1\dots b_s\hat b_{
s+1} K(b_{s+2}\dots b_r))$ corresponding to the $\psi(b_{s+1})$-term in $\Delta_\psi(\hat b_{s+1})$. Thus if $b_1\dots b_r$ is any monomial in $\pa_+ c$ then its contribution to the left hand side of \eqref{e:monomial1} is exactly
\[
K(b_1)\phi(b_2)\dots\phi(b_r)+\psi(b_1)K(b_2)\phi(b_3)\dots\phi(b_r)+\dots \psi(b_1)\dots\psi(b_{r-1})K(b_r)
\]
and \eqref{e:monomial1} follows.

We next consider monomials linear in $\hat\A_+$-generators contributing to \eqref{e:comp1} and \eqref{e:comp2}. To this end we define the maps $A,B,C,D\colon\A_+\to\A$ as follows on monomials.
\begin{align*}
A(b_1\dots b_s)&=\Gamma_\psi(\pa_+ b_1)K(b_2\dots b_s)\\
&+b_1\Gamma_\psi(\pa_+b_2)K(b_3\dots b_s)+\dots\\
&+b_1\dots b_{s-2}\Gamma_\psi(b_{s-1})K(b_s),
\end{align*}
if $s\ge 2$ and $0$ otherwise.
\begin{align*}
B(b_1\dots b_s)&=\hat b_1K(\pa_+(b_2\dots b_s))\\
&+b_1\hat b_2K(\pa_+(b_3\dots b_s))+\dots\\
&+b_1\dots b_{s-2}\hat b_{s-1}K(\pa_+(b_s)),
\end{align*}
if $s\ge 2$ and $0$ otherwise.
\begin{align*}
C(b_1\dots b_s)&=(\pa_+ b_1)\hat b_2K(b_3\dots b_s)\\
&+\pa_+(b_1b_2)\hat b_3K((b_4\dots b_s))+\dots\\
&+\pa(b_1\dots b_{s-2})\hat b_{s-1}K(b_s),
\end{align*}
if $s\ge 2$ and $0$ otherwise.
\begin{align*}
D(b_1\dots b_s)&=\hat b_1\pa_-(K(b_2\dots b_s))\\
&+b_1\hat b_2\pa_-(K(b_3\dots b_s))+\dots\\
&+b_1\dots b_{s-2}\hat b_{s-1}\pa_-(K(b_s)),
\end{align*}
if $s\ge 2$ and $0$ otherwise.

We have
\[
0=\Gamma_K(\pa_+\pa_+ c)=\Omega^{\Gamma_K\circ\pa_+}(\pa_+ c) + A(\pa_+ c) + B(\pa_+ c) + C(\pa_+ c).
\]
(To see this one subdivides contributing monomials as follows. Let $b_1\dots b_s$ be a monomial in $\pa_+c$. The first term corresponds to the $\hat \A_+$-variable and the $K$-variable both located in a $\pa b_j$-monomial for some $j$. The second term corresponds to the $\hat\A_+$-variable in a $\pa b_j$-monomial and the $K$-variable outside. The third term corresponds to the $\hat \A_+$-variable outside a $\pa_+ b_j$-monomial and the variable on which the $\pa_+$-operator acts being on the right of the $\hat \A_+$-generator. The fourth term corresponds to the $\hat \A_+$-variable outside a $\pa_+ b_j$-monomial and the variable on which the $\pa_+$-operator acts being on the left of the $\hat \A_+$-generator.)

Similarly, we have
\[
\Delta_\psi(\Gamma_K(\pa_+ c))_1= A(\pa_+ c) + C(\pa_+ c)+ D(\pa_+ c),
\]
where $w_1$ denotes the term of an element which is linear in the $\hat\A_+$-generators. Consequently, the contribution of monomials of the second kind to the sum of \eqref{e:comp1} and \eqref{e:comp2} is
\begin{align*}
&\Gamma_\phi(\pa_+ c)+\Gamma_\psi(\pa_+ c)+\Omega^{\Gamma_K\circ\pa_+}(\pa_+ c)+\Delta_\psi(\Gamma_K(\pa_+ c))_1\\
&=\Gamma_\phi(\pa_+ c)+\Gamma_\psi(\pa_+ c)+ B(\pa_+c) + D(\pa_+ c)=0.
\end{align*}
The lemma follows.
\end{pf}

\section{Examples}\label{S:examples}
The results of this paper allow us to construct many interesting Legendrian submanifolds. To illustrate this, we shall apply Theorem \ref{t:main} to some of the loops of Legendrian knots discussed in \cite{K}. In particular, we will concentrate on augmentations of the contact homology of the resulting Legendrian tori. Our examples give a proof of Theorem \ref{t:ex}.

In any situation when the Legendrian $L\subset J^{1}(\R)$, of Maslov class $r=0$, only has Reeb chords of non-negative grading (and $L$ is the base point of a closed loop that gives rise to the Legendrian torus $\Lambda$), the following observation applies. The ``long chords'' (hat variables) of $\Lambda$ have positive grading so they cannot be augmented. We have to select index $0$ chords $b$ of $L$ (``short chords'') on which the augmentation should take a non-zero value. Grading $1$ short chords impose the same relations on these as when searching for augmentations of $L$ itself. So augmentations of $\Lambda$ are in fact also augmentations of $L$. But there are also relations imposed by the grading $1$ long chords. Because of the absence of negative grading chords, these are of the form $\Delta\hat b=b+\mu(b)$, where $\mu$ is the monodromy of our loop. Hence, an augmentation of $L$ will be an augmentation of $\Lambda$ if and only if it is invariant under $\mu$.

In \cite{K}, a natural loop is described in the space of braid-positive Legendrian knots. The $(p,2)$ torus knots are its simplest special case, when conjugating a single crossing from one end of the braid to the other already results in a closed loop. We will examine the cases $p=3$ and $p=7$.

\begin{figure}
\labellist
\small%\hair 2pt
\pinlabel $a_2$ at 160 520
\pinlabel $a_2$ at 400 525
\pinlabel $a_2$ at 690 495
\pinlabel $a_2$ at 255 110
\pinlabel $a_2$ at 705 95
\pinlabel $a_1$ at 175 425
\pinlabel $a_1$ at 515 415
\pinlabel $a_1$ at 885 470
\pinlabel $a_1$ at 345 140
\pinlabel $a_1$ at 695 190
\pinlabel $b_1$ at 75 325
\pinlabel $b_2$ at 140 325
\pinlabel $b_3$ at 210 325
\pinlabel $b_1$ at 395 325
\pinlabel $b_2$ at 460 325
\pinlabel $b_3$ at 530 325
\pinlabel $b_1$ at 735 325
\pinlabel $b_2$ at 800 325
\pinlabel $b_3$ at 870 325
\pinlabel $b_1$ at 260 175
\pinlabel $b_2$ at 260 -5
\pinlabel $b_3$ at 330 -5
\pinlabel $b_2$ at 610 -5
\pinlabel $b_3$ at 675 -5
\pinlabel $c$ at 480 515
\pinlabel $c$ at 895 425
\pinlabel $c$ at 355 95
\pinlabel $c$ at 740 -5
\pinlabel $d$ at 555 450
\pinlabel $d$ at 785 440
\pinlabel $d$ at 150 145
\endlabellist
\centering
\includegraphics[width=.8\linewidth]{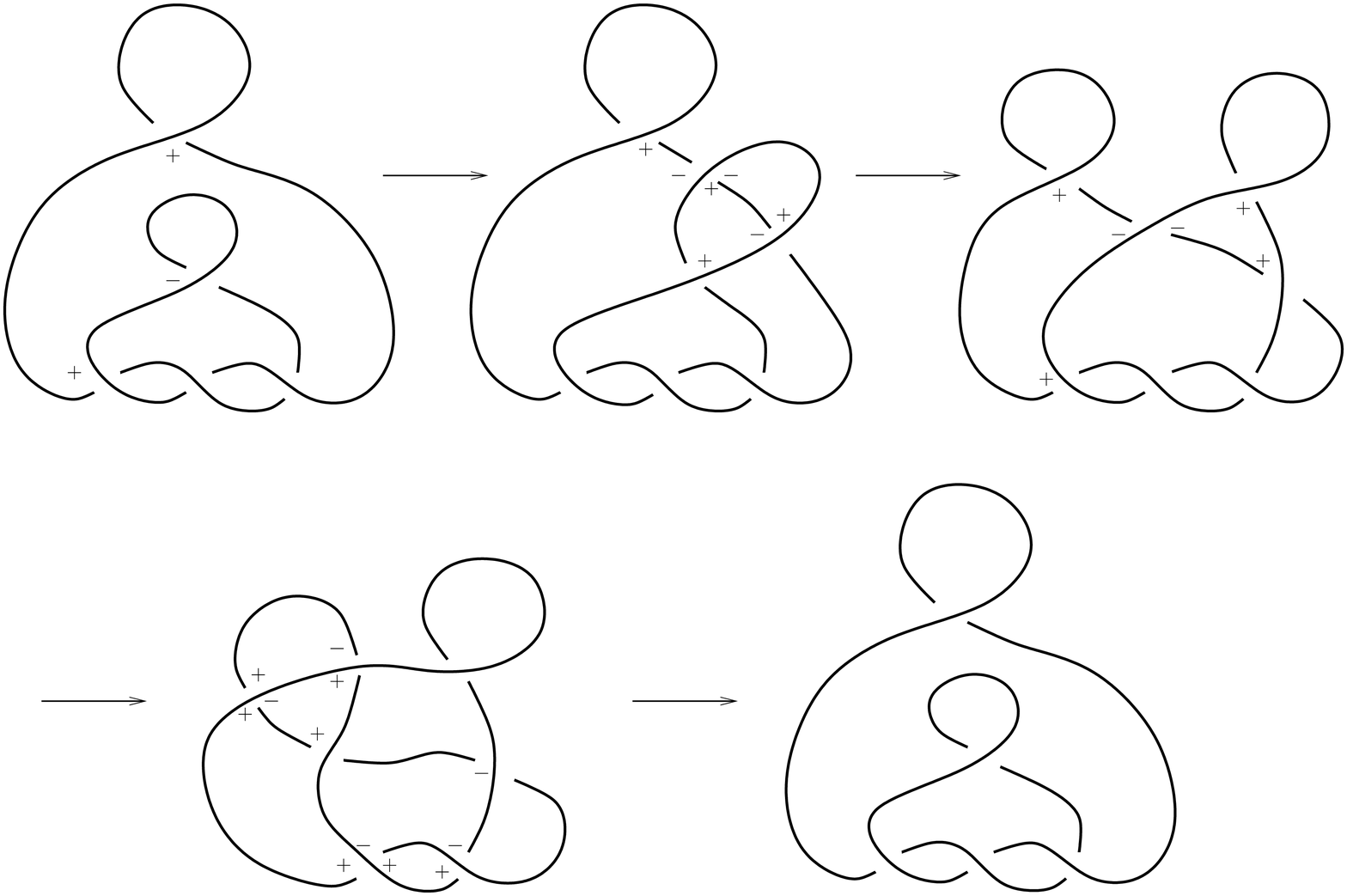}
\caption{A loop of Legendrian trefoil knots}
\label{fig:trefoilloop}
\end{figure}

\begin{pf}[Proof of Theorem \ref{t:ex} (a)]
The diagram of the loop for $p=3$ is reprinted in Figure \ref{fig:trefoilloop}. We recall from \cite[Section 5]{K} that, after restoring the original labels, the monodromy $\mu$ of this loop acts on the index $0$ variables $b_1$, $b_2$, and $b_3$ as follows:
\[\mu(b_1)=1+b_2b_3;\quad\mu(b_2)=b_1;\quad\mu(b_3)=b_2.\]
The crossings $a_1$ and $a_2$ are of index $1$ and they are not cycles, so for the time being we only concern ourselves with their boundaries (see, for example, \cite[Example 2.14]{K}):
\[\partial a_1=1+b_1+b_3+b_1b_2b_3;\quad\partial a_2=b_2+b_1b_2+b_2b_3+b_2b_3b_1b_2.\]
The other three relations that determine the index $0$ part $CH_0$ of the contact homology of the corresponding torus are
\begin{equation}\label{e:trefoil1}
\Delta\hat b_1=b_1+1+b_2b_3;\quad\Delta\hat b_2=b_2+b_1;\quad\Delta\hat b_3=b_3+b_2.
\end{equation}
The last two relations in \eqref{e:trefoil1} implies that $CH_0$ has a single generator. We denote it $b$. It is subject to the relations
\[1+b^3=0;\quad b+b^4=0;\quad 1+b+b^2=0.\]
The first two of these follows from the third by multiplication with $b+1$ and $b(b+1)$, respectively. Therefore,
\[CH_0\cong\Z_2[b]/\langle1+b+b^2=0\rangle.\]
Thus $CH_0$ is non-trivial. On the other hand because neither $0$ nor $1$ is a root of the relation, this contact homology cannot be augmented. (Indeed, the five augmentations of the trefoil are permuted in a single cycle by $\mu$, so none of them is fixed.)
\end{pf}

\begin{pf}[Proof of Theorem \ref{t:ex} (b)]
Let $p=7$ for the knots considered above. For convenience we introduce terminology that applies for any odd number $p$. The loop is the same as in Figure \ref{fig:trefoilloop} except that there are $p-3$ more index $0$ crossings $b_4,\ldots,b_p$ that, just like $b_2$ and $b_3$, are essentially unaffected by the monodromy (they simply get re-labeled at the end).

They do however influence $\partial a_1$ and $\partial a_2$ as follows:
\[\partial a_1=1+B_{11};\quad\partial a_2=1+B_{22}+B_{21}B_{12}.\]
Here, the polynomials $B_{ij}$ are natural generalizations of the expressions $B_{11}=b_1+b_3+b_1b_2b_3$, $B_{12}=1+b_2b_3$, $B_{21}=1+b_1b_2$, and $B_{22}=b_2$ in the $p=3$ case. See \cite[Section 6]{K} for more details. In particular when $p=7$, we obtain
\begin{eqnarray*}%\begin{split}
\Delta a_1=\partial a_1&=&1+b_1+b_3+b_5+b_7+b_1b_2b_3+b_1b_2b_5+b_1b_2b_7+b_1b_4b_5+b_1b_4b_7\\
&+&b_1b_6b_7+b_3b_4b_5+b_3b_4b_7+b_3b_6b_7+b_5b_6b_7+b_1b_2b_3b_4b_5+b_1b_2b_3b_4b_7\\
&+&b_1b_2b_3b_6b_7+b_1b_2b_5b_6b_7+b_1b_4b_5b_6b_7+b_3b_4b_5b_6b_7+b_1b_2b_3b_4b_5b_6b_7;\\
\Delta a_2=\partial a_2&=&1+b_2+b_4+b_6+b_2b_3b_4+b_2b_3b_6+b_2b_5b_6+b_4b_5b_6+b_2b_3b_4b_5b_6\\
&+&(1+b_2b_3+b_2b_5+b_2b_7+b_4b_5+b_4b_7+b_6b_7+b_2b_3b_4b_5\\
&&+b_2b_3b_4b_7+b_2b_3b_6b_7+b_2b_5b_6b_7+b_4b_5b_6b_7+b_2b_3b_4b_5b_6b_7)\\
&\cdot&(1+b_1b_2+b_1b_4+b_1b_6+b_3b_4+b_3b_6+b_5b_6+b_1b_2b_3b_4\\
&&+b_1b_2b_3b_6+b_1b_2b_5b_6+b_1b_4b_5b_6+b_3b_4b_5b_6+b_1b_2b_3b_4b_5b_6).
%\end{split}
\end{eqnarray*}
These formulas are special cases of  \cite[Theorem 6.7]{K}. The relations $\Delta \hat b_2=b_2+b_1=0,\ldots,\Delta\hat b_p=b_p+b_{p-1}=0$ reduce $CH_0$ of this torus, to a single-generator algebra. If $b$ denotes a generator of $CH_0$ then, when $p=7$, the previous two formulas impose the relations
\[1+b^7=0\text{ and }1+b+b^5+(1+b^4+b^6)^2=b+b^5+b^8+b^{12}=0.\]
The second of these is equal to $b+b^5$ times the first. For a general $p$, identifying the grading $0$ generators leads to the reductions
\[B_{11}\mapsto Q_p(b),\quad B_{12},B_{21}\mapsto Q_{p-1}(b),\quad B_{22}\mapsto Q_{p-2}(b),\]
where the polynomials $Q_k$ are defined by setting $Q_{-1}(b)=0$, $Q_0(b)=1$ and then applying the recursion $Q_{k}(b)=bQ_{k-1}(b)+Q_{k-2}(b)$. They are also characterized by the formula
\[\begin{bmatrix} b&1\\ 1&0\end{bmatrix}^k=\begin{bmatrix} Q_k(b)&Q_{k-1}(b)\\ Q_{k-1}(b)&Q_{k-2}(b)\end{bmatrix}.\]
(This follows by a straightforward induction argument, explicit formulas for the coefficients can be obtained from Pascal's triangle). In particular,
because the determinant of the above matrix is $1$, $\Delta a_2=\partial a_2=1+B_{22}+B_{21}B_{12}$ reduces to $1+Q_{p-2}(b)+Q_{p-1}(b)^2=Q_{p-2}(b)+Q_p(b)Q_{p-2}(b)=(1+Q_p(b))Q_{p-2}(b)$. Here of course, $1+Q_p(b)$ is just the reduction of $\Delta a_1=\partial a_1$. Thus so far, we obtained only one relation, $1+Q_p(b)=0$, for the single generator $b$ of $CH_0$.

The other relation, again for $p=7$, comes from
\begin{eqnarray*}
\Delta\hat b_1&=&b_1+1+b_2b_3+b_2b_5+b_2b_7+b_4b_5+b_4b_7+b_6b_7+b_2b_3b_4b_5\\
&+&b_2b_3b_4b_7+b_2b_3b_6b_7+b_2b_5b_6b_7+b_4b_5b_6b_7+b_2b_3b_4b_5b_6b_7,
\end{eqnarray*}
which simplifies to
\[1+b+b^4+b^6=0.\]
(For general $p$, $\Delta\hat b_1=b_1+B_{21}$, which reduces to $b+Q_{p-1}(b)$. This follows from Theorem \ref{t:main} and \cite[Proposition 8.2]{K}.)
The Euclidean algorithm shows that the greatest common divisor of this and $1+b^7$ is $1+b^2+b^3+b^4$. (Indeed, $(1+b+b^2)(1+b^2+b^3+b^4)=1+b+b^4+b^6$ and $(1+b^2+b^3)(1+b^2+b^3+b^4)=1+b^7$.) Thus in this case,
\begin{equation}\label{eq:7,2}
CH_0\cong\Z_2[b]/\langle1+b^2+b^3+b^4=0\rangle.
\end{equation}
Setting $b=1$ defines an augmentation, which corresponds to the fact that $\varepsilon(b_1)=\cdots=\varepsilon(b_7)=1$ is an invariant augmentation of the $(7,2)$ torus knot. (Other, more typical augmentations of the knot are not invariant under the monodromy of the loop: $81$ of them form nine $9$--cycles, and there is a $3$--cycle too.) With this augmentation of the torus, the linearized differential $\tilde\Delta$ takes the following form:
\begin{eqnarray}
&&\label{first}\tilde\Delta b_1=\cdots=\tilde\Delta b_7=0;\quad
\tilde\Delta a_1=b_1+b_4+b_7,\quad\tilde\Delta a_2=0;\\
&&\label{second}\tilde\Delta\hat b_1=b_1+b_3+b_6,\quad\tilde\Delta\hat b_2=b_1+b_2,\ldots,\tilde\Delta\hat b_7=b_6+b_7;\\
&&\label{third}\tilde\Delta\hat a_1=a_1+a_2+\hat b_1+\hat b_4+\hat b_7,\quad\tilde\Delta\hat a_2=0.
\end{eqnarray}

Here, \eqref{first} and \eqref{second} follow directly from the formulas already given, but \eqref{third} needs explanation. First, we need to understand $\mu(a_1)$ and $\mu(a_2)$. The key to this is \cite[Remark 3.4]{K} and the second diagram in Figure \ref{fig:trefoilloop}.
% put a + and - sign on the diagram
Note that there exist no admissible disks there with a positive corner at $a_1$ and a negative corner at the newly created $c$. Thus the image of $a_1$ is itself in this step, and it isn't affected later either, except for re-labeling to $\mu(a_1)=a_2$ at the end. There are, however, admissible disks from $a_2$ to $c$. With their contributions, the image of $a_2$ is $a_2+B_{21}d$. The two triangle moves that follow affect $d$, namely $d\mapsto d+a_1c\mapsto d+b_1a_2+a_1c$. Then at the last move before re-labeling, $d\mapsto 0$ and $b_1\mapsto B'_{21}$. Here, $B'_{21}$ refers to the braid with crossings labeled $b_2,\ldots,b_p,c$. So far, the image of $a_2$ is $a_2+B_{21}(B'_{21}a_2+a_1c)$. This gets re-labeled to $\mu(a_2)=a_1+B_{12}(B_{21}a_1+a_2b_p)$.

So we have
\[\Delta\hat a_1=a_1+a_2+\Gamma_\mu(\partial a_1)\]
and
\[\Delta\hat a_2=a_2+a_1+B_{12}B_{21}a_1+B_{12}a_2b_p+\Gamma_\mu(\partial a_2).\]
All monomials in these expressions have a single grading $1$ factor and several augmented grading $0$ factors. To linearize them, we just have to count the number of times each grading $1$ variable appears. When $p=7$, we have already checked that $B_{12}$ and $B_{21}$ are sums of an odd number of terms ($13$, to be exact). To compute the contributions from the $\Gamma_\mu(\partial a_i)$ terms, note that it is essentially the same task as the linearization of $\Delta a_i=\partial a_i$: each monomial contributes the sum of its terms to the latter, and the sum of the hat equivalents of its terms to the former.

Now it is straightforward to compute that the homology of this complex has rank $1$ in gradings $1$ and $2$ and rank $0$ everywhere else. This is identical to the linearized contact homology of the standard torus that is the trace of the constant isotopy of the unknot. (The latter has single generators in gradings $1$ and $2$ and a trivial differential, hence $0$ is its only augmentation.) Yet, the torus derived from the $(7,2)$ torus knot is different since its $CH_0$, given by \eqref{eq:7,2}, is non-trivial.
\end{pf}


\begin{thebibliography}{999}

\bibitem{B}
F. Bourgeois,
{\em A Morse-Bott approach to Contact Homology},
PhD thesis, Stanford University (2002)


\bibitem{BEHWZ}
F. Bourgeois, Y. Eliashberg, H. Hofer, K. Wysocki, E. Zehnder,
{\em Compactness results in symplectic field theory},
Geom. Topol. {\bf 7} (2003), 799--888


\bibitem{Ch}
Y. Chekanov,
{\em Differential algebra of Legendrian links},
Invent. Math. {\bf 150} (2002), no. 3, 441--483


\bibitem{E}
T. Ekholm,
{\em Morse flow trees and Legendrian contact homology in 1-jet spaces}, 
Geom. Topol. {\bf 11} (2007), 1083--1225

\bibitem{EES1}
T. Ekholm, J. Etnyre, M. Sullivan,
{\em Non-isotopic Legendrian submanifolds in ${\mathbb R}\sp {2n+1}$},
J. Differential Geom. {\bf 71} (2005), no. 1, 85--128

\bibitem{EES2}
T. Ekholm, J. Etnyre, M. Sullivan,
{\em The contact homology of Legendrian submanifolds in ${\mathbb R}\sp {2n+1}$}, J. Differential Geom. {\bf 71} (2005), no. 2, 177--305

\bibitem{EES3}
T. Ekholm, J. Etnyre, M. Sullivan,
{\em Orientations in Legendrian contact homology and exact Lagrangian immersions},
Internat. J. Math. {\bf 16} (2005), no. 5, 453--532

\bibitem{EES4}
T. Ekholm, J. Etnyre, M. Sullivan,
{\em Legendrian Contact Homology in $P\times\R$},
Trans. Amer. Math. Soc.  {\bf 359}  (2007),  no. 7, 3301--3335 

\bibitem{EHK} T. Ekholm, K. Honda, T. Kalman, 
{\em Legendrian knots and exact Lagrangian cobordisms}, 
in preparation

\bibitem{El}
Y. Eliashberg,
{\em Invariants in contact topology},
Proceedings of the International Congress of Mathematicians, Vol. II (Berlin, 1998). Doc. Math. 1998, Extra Vol. II, 327--338


\bibitem{EGH}
Y. Eliashberg, A. Givental, H. Hofer,
{\em Introduction to symplectic field theory},
GAFA 2000 (Tel Aviv, 1999). Geom. Funct. Anal. 2000, Special Volume, Part II, 560--673


\bibitem{FOOO}
K. Fukaya, Y.-G. Oh, H. Ohta, K. Ono,
{\em Lagrangian intersection Floer theory - anomaly and obstructon},
preprint (2000)


\bibitem{HWZI}
H. Hofer, K. Wysocki, E. Zehnder,
{\em A General Fredholm Theory I: A Splicing-Based Differential Geometry},
arXiv:math/0612604    

\bibitem{HWZII}
H. Hofer, K. Wysocki, E. Zehnder,
{\em A General Fredholm Theory II: Implicit Function Theorems},
arXiv:0705.1310

\bibitem{K}
T. K{\'a}lm{\'a}n, 
{\em Contact homology and one parameter families of Legendrian knots}, 
Geom. Topol. {\bf 9} (2005), 2013--2078


\bibitem{Ng1}
L. Ng, 
{\em Knot and braid invariants from contact homology. I.},
Geom. Topol. {\bf 9} (2005), 247--297

\bibitem{Ng2}
L. Ng,
{\em Knot and braid invariants from contact homology. II. With an appendix by the author and Siddhartha Gadgil},
Geom. Topol. {\bf 9} (2005), 1603--1637
\end{thebibliography}
\end{document}